\documentclass[12pt]{amsart}
\usepackage{amssymb}

\evensidemargin\paperwidth
\advance\evensidemargin-\textwidth
\advance\oddsidemargin-1in
\evensidemargin\oddsidemargin

\topmargin\paperheight
\advance\topmargin-\textheight
\advance\topmargin-1in

\title[Classification and Centralizers]{Concrete Classification and Centralizers of Certain $\mathbb{Z}^2 \rtimes {\rm SL}(2,\mathbb{Z})$-actions}
\author{Hiroki Sako}
\address{Department of Mathematical Sciences, University of Tokyo, Komaba, Tokyo, 153-8914, Japan}
\email{hiroki@ms.u-tokyo.ac.jp}

\newcommand{\C}{\mathbb{C}}
\newcommand{\R}{\mathbb{R}}
\newcommand{\Q}{\mathbb{Q}}
\newcommand{\Z}{\mathbb{Z}}

\newcommand{\tr}{\textrm{tr}}
\newcommand{\id}{{\rm id}}
\newcommand{\G}{\mathbb{Z}^2 \rtimes {\rm SL}(2,\mathbb{Z})}

\theoremstyle{plain}
\newtheorem{theorem}{Theorem}[section]
\theoremstyle{plain}
\newtheorem{proposition}[theorem]{Proposition}
\theoremstyle{plain}
\newtheorem{lemma}[theorem]{Lemma}
\theoremstyle{plain}
\newtheorem{corollary}[theorem]{Corollary}
\theoremstyle{plain}
\newtheorem{definition}[theorem]{Definition}
\theoremstyle{remark}
\newtheorem{remark}[theorem]{Remark}
\theoremstyle{remark}
\newtheorem{acknowledgment}{Acknowledgment}

\begin{document}
\begin{abstract}
We introduce a new class of actions of the group $\G$ on finite von Neumann algebras and call them twisted Bernoulli shift actions. We classify these actions up to conjugacy and give an explicit description of their centralizers. We also distinguish many of those actions on the AFD $\mathrm{II}_1$ factor in view of outer conjugacy.
\end{abstract}

\maketitle

\section{Introduction}
We consider the classification of $\G$-actions on finite von Neumann algebras in this paper. Mainly, we concentrate on the case that the finite von Neumann algebra is the AFD factor of type $\mathrm{II}_1$ or non-atomic abelian.

There are two difficulties for analyzing discrete group actions on operator algebras. The first is that we do not have various ways to construct actions. The second is that we can not analyze them by concrete calculation in most cases. To give many examples of actions which admit concrete analysis, we introduce a class of trace preserving $\G$-actions on finite von Neumann algebras and call them twisted Bernoulli shift actions. We classify those actions up to conjugacy and study them up to outer conjugacy.

An action $\beta(H, \mu, \chi)$ in the class is defined for a triplet $(H, \mu, \chi)$, where $H$ is an abelian countable discrete group, $\mu$ is a normalized scalar $2$-cocycle of $H$ and $\chi$ is a character of $H$. We obtain the action by restricting the so-called generalized Bernoulli shift action to a subalgebra $N(H, \mu)$ and ``twisting'' it by the character $\chi$. The process of restriction has a vital role in concrete analysis of these actions.

A $*$-isomorphism which gives conjugacy between two twisted Bernoulli shift actions $\beta(H_a, \mu_a, \chi_a)$ and $\beta(H_b, \mu_b, \chi_b)$ must be induced from an isomorphism between the two abelian groups $H_a$ and $H_b$. We prove this by concrete calculation (Section \ref{main theorems}). It turns out that there exist continuously many, non-conjugate $\G$-actions on the AFD factor of type $\mathrm{II}_1$ (Section \ref{examples}). By using the same technique, we describe the centralizers of all twisted Bernoulli shift actions. Here we should mention that the present work was motivated by the previous ones \cite{Ch}, \cite{NPS}, where similar studies were carried out in the case of ${\rm SL}(n, \Z)$.

In Section \ref{Section; outer conjugacy}, we distinguish many twisted Bernoulli shift actions in view of outer conjugacy. The classification for actions of discrete amenable groups on the AFD factor of type $\mathrm{II}_1$ was given by Ocneanu \cite{Oc}. Outer actions of countable amenable groups are outer conjugate. In the contrast to this, V. F. R. Jones \cite{Jo2} proved that any discrete non-amenable group has at least two non outer conjugate actions on the AFD factor of type $\mathrm{II}_1$. S. Popa (\cite{Po3}, \cite{Po4}, \cite{PoSa}, etc.) used the malleability/deformation arguments for the Bernoulli shift actions to study (weak) $1$-cocycles for the actions.
For some of twisted Bernoulli shift actions, which we introduce in this paper, it is shown that (weak) $1$-cocycles are represented in simple forms under some assumption on the (weak) $1$-cocycles. We prove that there exist continuously many twisted Bernoulli shift actions which are mutually non outer conjugate. This strengthens the above mentioned result due to Jones in the $\G$ cases.

\section{Preparations}
\subsection{Functions $\det$ and $\gcd$}
For the definition of twisted Bernoulli shift actions in Section \ref{definition}, we define two $\Z$-valued functions $\det$ and $\gcd$. The function $\det$ is given by the following equation:
\begin{eqnarray*}
   \det
   \left(
   \left(
   \begin{array}{c}
       q \\
       r
   \end{array}
   \right),
   \left(
   \begin{array}{c}
       q_0 \\
       r_0
   \end{array}
   \right)
   \right)
   = q r_0 - r q_0, \quad
   \left(
   \begin{array}{c}
       q \\
       r
   \end{array}
   \right),
   \left(
   \begin{array}{c}
       q_0 \\
       r_0
   \end{array}
   \right) \in \Z^2.
\end{eqnarray*}
The value of the function $\gcd$ at $k \in \Z^2$ is the greatest common divisor of the two entries. For $0 \in \Z^2$, let the value of gcd be $0$.

\begin{lemma}\label{det and gcd}\
\begin{enumerate}
\item
The action of ${\rm SL}(2,\Z)$ on $\Z^2$ preserves the functions $\det$ and $\gcd$, that is,
\begin{eqnarray*}
\det(k, k_0) &=& \det(\gamma \cdot k, \gamma \cdot k_0),\\
\gcd(k) &=& \gcd(\gamma \cdot k), \quad k, k_0 \in \Z^2, \gamma \in {\rm SL}(2, \Z).
\end{eqnarray*}
\item
The following equation holds true:
\begin{eqnarray*}
\det(k, k_0) &=& \gcd(k) + \gcd(k_0) - \gcd(k + k_0) \mod 2, \quad k, k_0 \in \Z^2.
\end{eqnarray*}
\end{enumerate}
\end{lemma}

\begin{proof}
The claim $(1)$ is a well-known fact, so we prove the claim $(2)$. For the function $\gcd$, we get
\begin{eqnarray*}
   \gcd
   \left(
     \left(
     \begin{array}{c}
        q \\
        r
     \end{array}
     \right)
   \right) =
   \left\{
   \begin{array}{l}
      1 \ \mod{2}, \quad ({\rm either \ } q \ {\rm or} \ r \ {\rm is \ odd}),\\
      0 \ \mod{2}, \quad ({\rm both \ } q \ {\rm and} \ r \ {\rm are \ even}).
   \end{array}
   \right.
\end{eqnarray*}
Since the action of ${\rm SL}(2, \Z)$ on $\Z^2$ preserves the functions $\det$ and $\gcd$, it suffices to show the desired equation against the following four pairs:
\begin{eqnarray*}
(k, k_0) &=&
     \left(
     \left(
     \begin{array}{c}
        0 \\
        0
     \end{array}
     \right),
     \left(
     \begin{array}{c}
        0 \\
        0
     \end{array}
     \right)
     \right),
     \left(
     \left(
     \begin{array}{c}
        1 \\
        0
     \end{array}
     \right),
     \left(
     \begin{array}{c}
        0 \\
        0
     \end{array}
     \right)
     \right),\\
   & &
     \left(
     \left(
     \begin{array}{c}
        1 \\
        0
     \end{array}
     \right),
     \left(
     \begin{array}{c}
        1 \\
        0
     \end{array}
     \right)
     \right),
     \left(
     \left(
     \begin{array}{c}
        1 \\
        0
     \end{array}
     \right),
     \left(
     \begin{array}{c}
        0 \\
        1
     \end{array}
     \right)
     \right) \ \mod{2}.
\end{eqnarray*}
\end{proof}

\subsection{ Scalar $2$-cocycles for abelian groups }\label{subsection; scalar 2-cocycles}
We fix some notations for countable abelian groups and their scalar $2$-cocycles. For the rest of this paper, let $H$ be an abelian countable discrete group and suppose that any scalar $2$-cocycle $\mu \colon H \times H \rightarrow \mathbb{T} =\{z \in \C \ | \ |z| = 1 \}$ is normalized, that is, $\mu(g, 0) = 1 = \mu(0, g)$ for $g \in H$.
We denote by $\mu^*$ the $2$-cocycle for $H$ given by
$\mu^*(g, h) = \overline{\mu(h, g)}, \ g, h \in H$.
Let $\mu^*\mu$ be the function on $H \times H$ defined by
\begin{eqnarray*}
   \mu^*\mu(g, h) = \overline{\mu(h, g)} \, \mu(g, h), \quad g, h \in H.
\end{eqnarray*}
This is a bi-character, that is, $\mu^*\mu(g, \cdot)$ and $\mu^*\mu(\cdot, h)$ are characters of $H$.
By using this function, we can describe the cohomology class of $\mu$. See \cite{OPT} for the proof of the following Proposition:

\begin{proposition}\label{bicharacter}
Two scalar $2$-cocycles $\mu_1$ and $\mu_2$ of $H$ are cohomologous if and only if $\mu_1^*\mu_1 = \mu_2^*\mu_2$.
\end{proposition}

Let $\C_\mu(H)$ be the twisted group algebra of $H$ with respect to the $2$-cocycle $\mu$. We denote by $\{ u_h \ |\  h \in H\}$ the standard basis for $\C_\mu(H)$ as $\C$-linear space. We recall that the $\C$-algebra $\C_\mu(H)$ has a structure of $*$-algebra defined by
\begin{eqnarray*}
   u_g \, u_h = \mu(g, h) \ u_{g + h}, \quad
   u_g^* = \overline{\mu(g, -g)} \, u_{-g}, \quad g, h \in H.
\end{eqnarray*}

Let $\widetilde\mu$ be the $\mathbb{T}$-valued function on $\oplus_{\Z^2} H \times \oplus_{\Z^2} H$ defined by
\begin{eqnarray}\label{equation; widetilde}
   \widetilde\mu(\lambda_1, \lambda_2) = \prod_{k \in \Z^2} \mu(\lambda_1(k), \lambda_2(k)), \quad \lambda_1, \lambda_2 \in \oplus_{\Z^2} H.
\end{eqnarray}
The function $\widetilde\mu$ is a normalized scalar $2$-cocycle for $\oplus_{\Z^2} H$.
Let $\Lambda(H)$ be the abelian group defined by
\begin{eqnarray*}
   \Lambda(H) = \left\{ \lambda : \Z^2 \rightarrow H \ \left| {\rm \ finitely \ supported \ and \ } \sum_{k \in \Z^2}\lambda(k) = 0 \right. \right\}.
\end{eqnarray*}
Its additive rule is defined by pointwise addition.

\subsection{Definition of a $\G$-action on $\Lambda(H)$}
The group ${\rm SL}(2, \Z)$ acts on $\Z^2$ as matrix-multiplication and the group $\Z^2$ also does on $\Z^2$ by addition. These two actions define the action of $\G$ on $\Z^2$ which is explicitly described as
\begin{eqnarray*}
   \left(
   \left(
   \begin{array}{c}
      q \\ r
   \end{array}
   \right)
   ,
   \left(
   \begin{array}{cc}
      x & y\\
      z & w
   \end{array}
   \right)
   \right)
\cdot
\begin{array}{c}
   \left(
   \begin{array}{c}
      q_0 \\ r_0
   \end{array}
   \right)
\end{array}
=
\left(
\begin{array}{c}
   q + x q_0 + y r_0 \\ r + z q_0 + w r_0
\end{array}
\right),
\end{eqnarray*}
for all
\begin{eqnarray*}
\left(
   \left(
   \begin{array}{c}
      q \\ r
   \end{array}
   \right)
   ,
   \left(
   \begin{array}{cc}
      x & y\\
      z & w
   \end{array}
   \right)
\right) \in \G, \quad
   \left(
   \begin{array}{c}
      q_0 \\ r_0
   \end{array}
   \right) \in \Z^2.
\end{eqnarray*}
We define an action of $\G$ on $\oplus_{\Z^2} H$ as
\begin{eqnarray*}
   (\gamma \cdot \lambda)(k) = \lambda(\gamma^{-1} \cdot k), \quad k \in \Z^2,
\end{eqnarray*}
for $\gamma \in \G$ and $\lambda \in \oplus_{\Z^2} H$.

\subsection{On the relative property (T) of Kazhdan}
We give the definition of the relative property (T) of Kazhdan for a pair of discrete groups.
\begin{definition}
Let $G \subset \Gamma$ be an inclusion of discrete groups. We say that the pair $(\Gamma, G)$ {\rm has the relative property (T)} if the following condition holds:\\
There exist a finite subset $F$ of $\Gamma$ and $\delta > 0$ such that if $\pi : \Gamma \rightarrow \mathcal{U(H)}$ is a unitary representation of $\Gamma$ on a Hilbert space $\mathcal{H}$ with a unit vector $\xi \in \mathcal{H}$ satisfying $ \| \pi(g)\xi - \xi \| < \delta$ for $g \in F$,
then there exists a non-zero vector $\eta \in \mathcal{H}$ such that $\pi(h)\eta =\eta$ for $h \in G$.
\end{definition}
Instead of this original definition, we use the following condition.
\begin{proposition}\label{Proposition; relative (T)}$($\cite{Joli}$)$
Let $G \subset \Gamma$ be an inclusion of discrete groups. The pair $(\Gamma, G)$ has the relative property (T) if and only if the following condition holds:\\
For any $\epsilon > 0$, there exist a finite subset $F$ of $\Gamma$ and $\delta > 0$ such that if $\pi : \Gamma \rightarrow \mathcal{U(H)}$ is a unitary representation of $\Gamma$ on a Hilbert space $\mathcal{H}$ with a unit vector $\xi \in \mathcal{H}$ satisfying $\| \pi(g)\xi - \xi \| < \delta$ for $g \in F$,
then $\| \pi(h)\xi - \xi \| < \epsilon$ for $h \in G$.
\end{proposition}
The pair $(\Z^2 \rtimes {\rm SL}(2,\Z), \Z^2)$ is a typical example of group with the relative property (T). See \cite{Bu} or \cite{Sh} for the proof.

\subsection{Weakly mixing actions}
An action of a countable discrete group $G$ on a von Neumann algebra $N$ is said to be \textit{ergodic} if any $G$-invariant element of $N$ is a scalar multiple of $1$. The {\it weak mixing property} is a stronger notion of ergodicity.
\begin{definition}\label{Definition; WM}
Let $N$ be a von Neumann algebra with a faithful normal state $\phi$. A state preserving action $(\rho_g)_{g \in G}$ of a countable discrete group $G$ on $N$ is said to be {\rm weakly mixing} if for every finite subset $\{a_1, a_2, \hdots, a_n\} \subset N$ and $\epsilon > 0$, there exists $g \in G$ such that $| \phi(a_i\rho_g(a_j)) - \phi(a_i) \phi(a_j)| < \epsilon, \  i, j = 1, \hdots ,n$.
\end{definition}
The following is a basic characterization of the weak mixing property.  Between two von Neumann algebra $N$ and $M$, $N \otimes M$ stands for the tensor product von Neumann algebra.
\begin{proposition}[Proposition D.2 in \cite{Vaes}]\label{Weakly Mixing}
Let a countable discrete group $G$ act on a finite von Neumann algebra $(N, \tr)$ by trace preserving automorphisms $(\rho_g)_{g \in G}$. The following statements are equivalent:
\begin{enumerate}
   \item The action $(\rho_g)$ is weakly mixing.
   \item The only finite-dimensional invariant subspace of $N$ is $\C1$.
   \item For any action $(\alpha_g)$ of $G$ on a finite von Neumann algebra $(M, \tau)$, we have $(N \otimes M)^{\rho \otimes \alpha} = 1 \otimes M^\alpha$, where $(N \otimes M)^{\rho \otimes \alpha}$ and $M^\alpha$ are the fixed point subalgebras.
\end{enumerate}
\end{proposition}

\subsection{A remark on group von Neumann algebras}\label{Algebraical computation}
Let $\Gamma$ be a discrete group and let $\mu$ be a scalar $2$-cocycle of a countable group $\Gamma$. A group $\Gamma$ acts on the Hilbert space $\ell^2 \Gamma$ by the following two ways;
\begin{eqnarray*}
	u_\gamma(\delta_g) = \mu(\gamma, g) \delta_{\gamma g}, \quad \rho_\gamma(\delta_g) = \mu(g, \gamma^{-1} ) \delta_{g \gamma^{-1}}, \quad \gamma, g \in \Gamma.
\end{eqnarray*}
These two representations commute with each other. The von Neumann algebra $L_\mu(\Gamma)$ generated by the image of $u$ is called the group von Neumann algebra of $\Gamma$ twisted by $\mu$. The normal state $\langle \cdot \delta_e, \delta_e \rangle$ is a trace on $L_\mu(\Gamma)$. The vector $\delta_e$ is separating for $L_\mu(\Gamma)$. For any element $a \in L_{\mu}(\Gamma)$, we define the square summable function $a(\cdot)$ on $\Gamma$ by $a \delta_e = \sum a(g) \delta_g$. The function $a(\cdot)$ is called the {\it Fourier coefficient} of $a$.
We write $a = \sum_{g \in \Gamma} a(g) u_g$ and call this the Fourier expansion of $a$.
The Fourier expansion of $a^*$ is given by $a^* = \sum_{g \in \Gamma} \overline{\mu(g, g^{-1}) a(g^{-1})} u_g$,
since the Fourier coefficient $a^*(g) = \langle a^* \delta_e, \delta_g \rangle$ is described as
\begin{eqnarray*}
	\left\langle \delta_e, a \rho_{g^{-1}} \delta_e \right\rangle
	=  \left\langle \rho_{g^{-1}}^* \delta_e, \sum_g a(g) \delta_g \right\rangle
	=  \overline{\mu(g, g^{-1}) a(g^{-1})}.
\end{eqnarray*}
Here we used the equation $\rho_{g^{-1}}^* = \overline{\mu(g, g^{-1})} \rho_g$, which is verified by direct computation. For two elements $a,b$, the Fourier coefficient of $ab$ is given by
\begin{eqnarray*}
	a b(\gamma)  = \langle b \delta_e, \rho_{\gamma^{-1}} a^* \delta_e \rangle
	             = \sum_g \overline{a^*(g) \mu(g, \gamma)} b(g \gamma)
                 = \sum_g \mu(g^{-1}, g \gamma) a(g^{-1}) b(g \gamma).
\end{eqnarray*}
This equation allows us to calculate the Fourier coefficient algebraically, that is,
\begin{eqnarray*}
    a b = \sum_{\gamma} \left( \sum_g \mu(g^{-1}, g \gamma) a(g^{-1}) b(g \gamma) \right) u_\gamma
        = \sum_{\gamma} \left( \sum_{g h = \gamma} a(g) b(h) \right) u_\gamma.
\end{eqnarray*}

For a subgroup $\Lambda \subset \Gamma$, the subalgebra $\{ u_\lambda \ | \ \lambda \in \Lambda \}^{\prime\prime} \subset L_\mu(\Gamma)$ is isomorphic to $L_\mu(\Lambda)$. We sometimes identify them. An element $a \in L_\mu(\Gamma)$ is in the subalgebra $L_\mu(\Lambda)$ if and only if the Fourier expansion $a(\cdot) \colon \Gamma \rightarrow \C$ is supported on $\Lambda$, since the trace $\langle \cdot \delta_e, \delta_e \rangle$ preserving conditional expectation $E$ from $L_\mu(\Gamma)$ onto $L_\mu(\Lambda)$ is described as $E(a) = \sum_{\lambda \in \Lambda} a(\lambda) u_\lambda$.

\section{ Definition of twisted Bernoulli shift actions}\label{definition}
In this section, we introduce twisted Bernoulli shift actions of $\G$ on finite von Neumann algebras. The action is defined for a triplet $i = (H, \mu, \chi)$, where $H \neq \{ 0 \}$ is an abelian countable discrete group, $\mu$ is a normalized scalar $2$-cocycle of $H$ and $\chi$ is a character of $H$. The finite von Neumann algebra, on which the group $\G$ acts, is defined by the pair $(H, \mu)$.

We introduce a group structure on the set $\Gamma_0 = \widehat{H} \times \Z^2 \times {\rm SL}(2, \Z)$ as
\begin{eqnarray*}
    (c_1, k, \gamma_1)(c_2, l, \gamma_2) = \left(c_1 c_2 \chi^{\det(k, \gamma_1 \cdot l)}, k + \gamma_1 \cdot l, \gamma_1 \gamma_2 \right)
\end{eqnarray*}
for any $c_1, c_2 \in \widehat{H} , k, l \in \Z^2, \gamma_1, \gamma_2 \in {\rm SL}(2, \Z)$. The associativity is verified by Lemma \ref{det and gcd}. It turns out that the subsets $\widehat{H} = \widehat{H} \times \{ 0 \} \times \{ e \}$ and $G_0 = \widehat{H} \times \Z^2 \times \{ e \}$ are subgroups in $\Gamma_0$. It is easy to see that $G_0$ is a normal subgroup of $\Gamma_0$ and that $\widehat{H}$ is a normal subgroup of $G_0$ and $\Gamma_0$. We get a normal inclusion of groups $G_0 / \widehat{H} \subset \Gamma_0 / \widehat{H}$ and this is isomorphic to $\Z^2 \subset \G$.

Before stating the definition of the twisted Bernoulli shift action, we define a $\Gamma_0$-action $\rho$ on the von Neumann algebra $L_{\widetilde{\mu}}(\oplus_{\Z^2} H)$. We denote by $u(\lambda) \in L_{\widetilde{\mu}}(\oplus_{\Z^2} H)$ the unitary corresponding to $\lambda \in \oplus_{\Z^2} H$. We define a faithful normal trace $\tr$ of $L_{\widetilde{\mu}}(\oplus_{\Z^2} H)$ in the usual way. For $c \in \widehat{H}, k \in \Z^2, \gamma \in {\rm SL}(2, \Z)$, let $\rho(c), \rho(k), \rho(\gamma)$ be the linear transformations on $\C_{\widetilde{\mu}}(\oplus_{\Z^2} H)$ given by,
\begin{eqnarray*}
   \rho(c)(u(\lambda)) &=& \left( \prod_{l \in \Z^2}c(\lambda(l)) \right) u(\lambda),\\
   \rho(k)(u(\lambda)) &=& \left( \prod_{m \in \Z^2}\chi(\lambda(m))^{\det(k, m)} \right) u(k \cdot \lambda),\\
   \rho(\gamma)(u(\lambda)) &=& u(\gamma \cdot \lambda), \quad \lambda \in \Lambda(H), \end{eqnarray*}
These maps are compatible with the multiplication rule and the $*$-operation of $\C_{\widetilde{\mu}}(\oplus_{\Z^2} H)$ . Since these maps preserve the trace, they extend to $*$-automorphisms on $L_{\widetilde{\mu}}(\oplus_{\Z^2} H)$. It is immediate to see that $\rho(c)$ commutes with $\rho(k)$ and $\rho(\gamma)$. For $k, l \in \Z^2$, we have the following relation:
\begin{eqnarray*}
    \rho(k) \circ \rho(l)(u(\lambda))
    &=& \prod_{m \in \Z^2} \chi(\lambda(m))^{\det(l, m)} \rho(k) (u(l \cdot \lambda))\\
    &=& \prod_{m \in \Z^2} \chi(\lambda(m))^{\det(l, m)}
               \prod_{m \in \Z^2} \chi((l \cdot \lambda)(m))^{\det(k, m)}
               u(k \cdot (l \cdot \lambda))\\
    &=& \prod_{m \in \Z^2} \chi(\lambda(m))^{\det(l, m)} \chi(\lambda(m))^{\det(k, m + l)}
    u((k + l) \cdot \lambda).
\end{eqnarray*}
By $\det(l, m) + \det(k, m + l) = \det(k, l) + \det(k + l, m)$, this equals to
\begin{eqnarray*}
    & & \left( \prod_{m \in \Z^2} \chi(\lambda(m)) \right)^{\det(k, l)}
               \prod_{m \in \Z^2} \chi(\lambda(m))^{\det(k + l, m)}
               u((k + l) \cdot \lambda)\\
    &=& \rho(\chi^{\det(k, l)}) \circ \rho(k + l)(u(\lambda)).
\end{eqnarray*}
Since $\det$ is $\mathrm{SL}(2, \Z)$-invariant (Lemma \ref{det and gcd}), for $k \in \Z^2, \gamma \in {\rm SL}(2, \Z)$, we get
\begin{eqnarray*}
      \rho(\gamma \cdot k) \circ \rho(\gamma)(u(\lambda))
    &=& \rho(\gamma \cdot k)(u(\gamma \cdot \lambda))\\
   &=& \prod_{l \in \Z^2} \chi((\gamma \cdot \lambda)(l))^{\det(\gamma \cdot k, l)}
            u((\gamma \cdot k) \cdot (\gamma \cdot \lambda))\\
   &=&  \prod_{l \in \Z^2} \chi(\lambda(l))^{\det(\gamma \cdot k, \gamma \cdot l)}
            u(\gamma \cdot (k \cdot \lambda))\\
   &=& \prod_{l \in \Z^2} \chi(\lambda(l))^{\det(k, l)} \rho(\gamma)(u(k \cdot \lambda))\\
   &=& \rho(\gamma) \circ \rho(k) (u(\lambda)), \quad \lambda \in \Lambda(H).
\end{eqnarray*}
By using the above two equations, $\rho$ satisfies the following formula:
\begin{eqnarray*}
& &(\rho(c_1) \circ \rho(k) \circ \rho(\gamma_1)) \circ (\rho(c_2) \circ \rho(l) \circ \rho(\gamma_2))\\
&=& \rho(c_1) \circ \rho(c_2) \circ \rho(k) \circ \rho(\gamma_1) \circ \rho(l) \circ \rho(\gamma_2)\\
&=& \rho(c_1) \circ \rho(c_2) \circ \rho(k) \circ \rho(\gamma_1 \cdot l) \circ \rho(\gamma_1) \circ \rho(\gamma_2)\\
&=& \rho(c_1) \circ \rho(c_2) \circ \rho(\chi^{\det(k, \gamma_1 \cdot l)}) \circ \rho(k + (\gamma_1 \cdot l)) \circ \rho(\gamma_1 \gamma_2)\\
&=& \rho(c_1 c_2 \chi^{\det(k, \gamma_1 \cdot l)}) \circ \rho(k + (\gamma_1 \cdot l)) \circ \rho(\gamma_1 \gamma_2).
\end{eqnarray*}
With $\rho(c, k, \gamma) = \rho(c) \circ \rho(k) \circ \rho(\gamma)$, $\rho$ gives a $\Gamma_0$-action on $L_{\widetilde{\mu}}(\oplus_{\Z^2} H)$.

We define the finite von Neumann algebra $N(H, \mu)$ as the group von Neumann algebra $L_{\widetilde{\mu}}(\Lambda(H))$. By using Fourier coefficients, we can prove that $N(H, \mu) $ is the fixed point algebra under the $\widehat{H}$-action $\rho(\widehat{H}, 0, e)$ on $L_{\widetilde{\mu}}(\oplus_{\Z^2} H)$. We get a $\G$-action on $N(H, \mu)$ by
\begin{eqnarray*}
    \beta(k, \gamma)(x) = \rho(1, k, \gamma)(x), \quad k \in \Z^2, \gamma \in {\rm SL}(2, \Z), x \in N(H, \mu).
\end{eqnarray*}
This is the definition of the twisted Bernoulli shift action $\beta = \beta(H, \mu, \chi)$ on $N(H, \mu)$.

We obtained the actions $\beta(H, \mu, \chi)$ not only by twisting generalized Bernoulli shift actions but also restricting to subalgebras $N(H, \mu) \subset L_{\widehat\mu}(\oplus_{\Z^2} H) = \bigotimes_{\Z^2} L_\mu(H)$.
This restriction allows us to classify the actions up to conjugacy in the next section. In order to give a variety of the actions, we twisted the shift actions by the character $\chi$ of the abelian group $H$.

\begin{remark}\label{remark; weak mixing}
The action $\beta | _{\Z^2} = \beta(H, \mu, \chi)| _{\Z^2}$ has the weak mixing property.
In definition \ref{Definition; WM}, we may assume that the Fourier coefficients of $a_i\ (i = 1,2, \cdots, n)$ are finitely supported, by approximating in the $L^2$-norm. Then for appropriate $k \in \Z^2$, we get $\tr(a_i \beta(k)(a_j)) = \tr(a_i)\tr(a_j), \ i,j =1,2,\cdots,n$.
\end{remark}

\section{Classification up to conjugacy}\label{main theorems}
In this section, we classify the twisted Bernoulli shift actions $ \{ \beta(H, \mu, \chi) \} $ up to conjugacy (Theorem \ref{thm1}). We prove that an isomorphism which gives conjugacy between two twisted Bernoulli shift actions is of a very special form. In fact it comes from an isomorphism in the level of base groups $H$. We also determine the centralizer of the $\G$-action $\beta(H, \mu, \chi)$ on $N(H,\mu)$ (Theorem \ref{thm2}).

We fix some notations for the proofs. We define $0, e_1, e_2 \in \Z^2$ as
\begin{eqnarray*}
   0 =
   \left(
   \begin{array}{c}
      0 \\ 0
   \end{array}
   \right), \
   e_1 =
   \left(
   \begin{array}{c}
      1 \\ 0
   \end{array}
   \right), \
   e_2 =
   \left(
   \begin{array}{c}
      0 \\ 1
   \end{array}
   \right).
\end{eqnarray*}
Let $\xi$ be the element of $\G$ satisfying
\begin{eqnarray*}
   \xi \cdot 0 = e_1, \quad
   \xi \cdot e_1 = e_2, \quad
   \xi \cdot e_2 = 0.
\end{eqnarray*}
The elements $\xi$ and $\xi^2$ are explicitly described as
\begin{eqnarray*}
   \left\{
   \begin{array}{lcl}
   \xi &=&
   \left(
   \begin{array}{cc}
      e_1,&
             \left(
             \begin{array}{rr}
                -1 & -1\\
                 1 &  0
             \end{array}
             \right)
   \end{array}
   \right),  \\[3mm]
   \xi^2 &=&
   \left(
   \begin{array}{cc}
      e_2,&
             \left(
             \begin{array}{rr}
                 0 &  1\\
                -1 & -1
             \end{array}
             \right)
   \end{array}
   \right).
   \end{array}
   \right.
\end{eqnarray*}
The order of $\xi$ is 3. Let $\eta, \delta \in {\rm SL}(2, \Z)$ be given by $
\eta =
\left(
\begin{array}{cc}
   -1& 0\\
    0&-1
\end{array}
\right)$,
$
\delta =
\left(
\begin{array}{cc}
   1&1\\
   0&1
\end{array}
\right)
$.
Let $D$ be the subset of all elements of $\Z^2$ fixed under the action of $\delta$, that is,
\begin{eqnarray*}
D =
\left\{
\begin{array}{c|c}
   \left(
   \begin{array}{c}
      n\\
      0
   \end{array}\right)
   &
   n \in \Z
\end{array}
\right\}.
\end{eqnarray*}
Then we get
\begin{eqnarray*}
\begin{array}{cc}
\begin{array}{llc}
\xi \cdot D =
\left\{
\begin{array}{c|c}
   \left(
   \begin{array}{c}
      1 - n\\
      n
   \end{array}\right)
   &
   n \in \Z
\end{array}
\right\}, &
\xi^2 \cdot D =
\left\{
\begin{array}{c|c}
   \left(
   \begin{array}{c}
      0\\
      1 - n
   \end{array}\right)
   &
   n \in \Z
\end{array}
\right\}.
\end{array}
\end{array}
\end{eqnarray*}
We define the subgroup $\Lambda_D(H)$ of $\Lambda(H)$ by
\begin{eqnarray*}
   \Lambda_D(H) = \{ \lambda \in \Lambda(H) \ | \ \lambda : \Z^2 \rightarrow H \ {\rm is \ supported \ on \ } D \}.
\end{eqnarray*}

Let $(H_a, \mu_a, \chi_a)$ and $(H_b, \mu_b, \chi_b)$ be triplets of countable abelian groups, their normalized $2$-cocycles and characters. For $h \in H_a$, we define $\lambda_h \in \Lambda_D(H_a)$ as
\begin{eqnarray*}
\lambda_h(k) =
\left\{
\begin{array}{rl}
   h & (k = e_1), \\
  -h & (k = 0), \\
   0 & (k \neq e_1, 0).
\end{array}
\right.
\end{eqnarray*}
For $g \in H_b$, we define $\sigma_g \in \Lambda_D(H_b)$ as
\begin{eqnarray*}
\sigma_g(k) =
\left\{
\begin{array}{rl}
   g & (k = e_1), \\
  -g & (k = 0), \\
   0 & (k \neq e_1, 0).
\end{array}
\right.
\end{eqnarray*}
We denote by $v(\sigma) \in N(H_b, \mu_b)$ the unitary corresponding to $\sigma \in \Lambda(H_b)$.

\begin{theorem}\label{thm1}
If $\pi \colon N(H_a, \mu_a) \rightarrow N(H_b, \mu_b)$ is a $*$-isomorphism giving conjugacy between $\beta_a = \beta(H_a, \mu_a, \chi_a)$ and $\beta_b = \beta(H_b, \mu_b, \chi_b)$, then there exists a group isomorphism $\phi = \phi_\pi \colon H_a \rightarrow H_b$ satisfying
\begin{enumerate}
\item\label{condition; phi circ}
$\pi(u(\lambda)) = v(\phi \circ \lambda)\ \mathrm{mod} \ \mathbb{T}$ for $\lambda \in \Lambda(H_a)$,
\item\label{condition; cohomologous}
the $2$-cocycles $\mu_a(\cdot, \cdot)$ and $\mu_b(\phi(\cdot), \phi(\cdot))$ of $H_a$ are cohomologous,
\item\label{condition; character}
$\chi_a^2 = (\chi_b \circ \phi)^2.$
\end{enumerate}
Conversely, given a group isomorphism $\phi \colon H_a \rightarrow H_b$ satisfying $(\ref{condition; cohomologous})$ and $(\ref{condition; character})$, there exists a $*$-isomorphism $\pi = \pi_\phi \colon N(H_a, \mu_a) \rightarrow N(H_b, \mu_b)$ which satisfies condition $(\ref{condition; phi circ})$ and gives conjugacy between $\beta_a$, $\beta_b$.
\end{theorem}

We note that by Proposition \ref{bicharacter} condition $(\ref{condition; cohomologous})$ for $\phi$ is equivalent to
\begin{enumerate}
\item[$(2)'$]\label{condition; cohomologous prime}
$\mu_a^* \mu_a(g, h) = \mu_b^* \mu_b(\phi(g), \phi(h)), \ g, h \in H_a.$
\end{enumerate}

\begin{proof}[Proof for the first half of Theorem \ref{thm1}]\ \\
Suppose that there exists a (not necessarily trace preserving) $*$-isomorphism $\pi$ from $N(H_a, \mu_a)$ onto $N(H_b, \mu_b)$ such that $\pi \circ \beta_a(\gamma) = \beta_b(\gamma) \circ \pi, \ \gamma \in \G$.

We prove that for every $h \in H_a$ there exists $\phi(h) \in H_b$ satisfying $\pi(u(\lambda_h)) = v(\sigma_{\phi(h)}) \mod{\mathbb{T}}$. Let $U_h$ denote the unitary in $N(H_b, \mu_b)$
\begin{eqnarray*}
   U_h = \pi\left(\overline{\mu_a(h, -h)} \, u(\lambda_h) \right), \quad h \in H_a.
\end{eqnarray*}
We identify $N(H, \mu)$ with the subalgebra of the infinite tensor product $\bigotimes_{\Z^2} L_\mu (H)$, which is canonically isomorphic to $L_{\widetilde\mu} \left( \oplus_{\Z^2} H \right)$.
The preimage $\pi^{-1}(U_h)$ can be written as $u_h^* \otimes u_h$. Here $u_h$ is the unitary corresponding to $h \in H_a$ and placed on $1 \in \Z^2$ and the unitary $u_h^*$ is placed on $0 \in \Z^2$. We describe $U_h$ as the Fourier expansion
$U_h = \sum_{\sigma \in \Lambda(H_b)} c(\sigma) v(\sigma)$.
Since $e_1$ and $0$ are fixed under the action of $\delta$, one has
\begin{eqnarray*}
\beta_b(\delta)^n(U_h) = \pi \circ \beta_a(\delta)^n (\pi^{-1} (U_h)) = U_h.
\end{eqnarray*}
It follows that the Fourier expansion $U_h = \sum_{\sigma \in \Lambda(H_b)} c(\sigma) v(\sigma)$ must satisfy that $c(\sigma) = c(\delta^{-n} \cdot \sigma)$ for every $\sigma \in \Lambda(H_b)$
and $n \in \Z$. For $\sigma \in \Lambda(H_b) \setminus \Lambda_D(H_b)$, the orbit of $\sigma$ under the action of $\delta^{-1}$ is an infinite set, since the support ${\rm supp}(\sigma) \subset \Z^2$ is not included in $D$.
It turns out that $c(\sigma) = 0$ for all $\sigma \in \Lambda(H_b) \setminus \Lambda_D(H_b)$ due to $\Sigma |c(\sigma)|^2 = 1 < +\infty$, so that
$U_h = \sum_{\sigma \in \Lambda_D(H_b)} c(\sigma) v(\sigma)$.

The unitary $\chi_a(h) U_h^*$ is also fixed under the action of $\delta$ and can be written as
\begin{eqnarray*}
    \chi_a(h) U_h^*
&=& \pi \left(\chi_a(h) \overline{\mu_a(h, -h)} u(-\lambda_h)\right)\\
&=& \pi \left(\chi_a(h) \overline{\mu_a(h, -h)} u(\xi \cdot \lambda_h)\right) \,
                       \pi\left(\overline{\mu_a(h, -h)} u(\xi^2 \cdot \lambda_h)\right)\\
&=& \beta_b(\xi)(U_h) \, \beta_b(\xi^2)(U_h).
\end{eqnarray*}
Letting $n e_1 = (n, 0)^T \in \Z^2$, we get
\begin{eqnarray*}
\beta_b(\xi)(U_h)
&=& \beta_b(\xi) \left( \sum c(\sigma) \, v(\sigma) \right)
= \sum_{\sigma \in \Lambda_D(H_b)} c(\sigma) \, v(\xi \cdot \sigma)
                                    \prod_{n \in \Z} \chi_b(\sigma(n e_1))^n,\\
\beta_b(\xi^2)(U_h)
&=& \beta_b(\xi^2) \left( \sum c(\sigma) \, v(\sigma) \right)
= \sum_{\sigma \in \Lambda_D(H_b)} c(\sigma) \, v(\xi^2 \cdot \sigma).
\end{eqnarray*}
Since Fourier expansion admits algebraical calculation as in subsection \ref{Algebraical computation}, the expansion of $\chi_a(h) U_h^*$ is
\begin{eqnarray*}
& & \chi_a(h) U_h^* = \beta_b(\xi)(U_h) \, \beta_b(\xi^2)(U_h)\\
&=& \sum_{\sigma_1, \sigma_2 \in \Lambda_D(H_b)} c(\sigma_1) \, c(\sigma_2) \, v(\xi \cdot \sigma_1) \, v(\xi^2 \cdot \sigma_2) \prod_{n \in \Z} \chi_b(\sigma_1(n e_1))^n \\
&=& \sum_{\sigma_1, \sigma_2 \in \Lambda_D(H_b)} c(\sigma_1) \, c(\sigma_2) \, \widetilde\mu_b(\xi \cdot \sigma_1, \xi^2 \cdot \sigma_2) \, \prod_{n \in \Z} \chi_b(\sigma_1(n e_1))^n \, v(\xi \cdot \sigma_1 + \xi^2 \cdot \sigma_2).
\end{eqnarray*}
The map
$
\Lambda_D(H_b) \times \Lambda_D(H_b) \ni (\sigma_1, \sigma_2) \mapsto \xi \cdot \sigma_1 + \xi^2 \cdot \sigma_2 \in \Lambda(H_b)
$
is injective. Indeed, $\sigma_1$ is uniquely determined by $\xi \cdot \sigma_1 + \xi^2 \cdot \sigma_2$, since $\sigma_1(k) = ( \xi \cdot \sigma_1 + \xi^2 \cdot \sigma_2 ) (\xi \cdot k), k \in D \setminus \{e_1\}$ and $\sigma_1(e_1) = - \sum_{k \in D \setminus \{e_1\}} \sigma_1(k)$. Here we used the condition $\sum \sigma_1(k) = 0$. The element $\sigma_2$ is also determined by $\xi \cdot \sigma_1 + \xi^2 \cdot \sigma_2$. Thus the index $(\sigma_1, \sigma_2)$ uniquely determines $\xi \cdot \sigma_1 + \xi^2 \cdot \sigma_2$.

We take arbitrary elements $\sigma_1, \sigma_2 \in \Lambda_D(H_b)$ and suppose that $c(\sigma_1) \neq 0, c(\sigma_2) \neq 0$. Since the unitary $\chi_a(h) U_h^*$ is invariant under the action of $\delta$ and the coefficient of $\xi \cdot \sigma_1 + \xi^2 \cdot \sigma_2$ is not zero, $\xi \cdot \sigma_1 + \xi^2 \cdot \sigma_2$ is supported on $D$. It follows that the elements $\sigma_1$ and $\sigma_2$ can be written as $\sigma_1 = \sigma_{\phi(h)} = \sigma_2$, by some $\phi(h) \in H_b$. Indeed, since the subsets $D \setminus \{0, e_1\}, \xi D \setminus \{e_1, e_2\}$ and $\xi^2 D \setminus \{e_2, 0\}$ are mutually disjoint, the element $\xi \cdot \sigma_1$ must be supported on $\{e_1, e_2\}$ and the element $\xi^2 \cdot \sigma_2$ must be supported on $\{e_2, 0\}$. By the assumption $\sum_{k \in \Z^2} \sigma_i(k) = 0\ (i = 1, 2)$, $\sigma_i$ can be written as $\sigma_{\phi(h_i)}$. Then using the fact that $(\xi \cdot \sigma_1 + \xi^2 \cdot \sigma_2)(e_2) = \sigma_1 (e_1) + \sigma_2(0) = 0$, we get that $\sigma_1 = \sigma_{\phi(h)} = \sigma_2$ for some $h \in H_b$.
This means that there exists only one $\sigma \in \Lambda(H_b)$ such that $c(\sigma) \neq 0$ and that it is of the form $\sigma = \sigma_{\phi(h)}$. Then the unitary $U_h$ satisfies $U_h = \pi(u(\lambda_h)) = v(\sigma_{\phi(h)}) \mod{\mathbb{T}}$.

We claim that the map $\phi = \phi_\pi : H_a \rightarrow H_b$ is a group isomorphism. For all $h_1, h_2 \in H_a$, we get
\begin{eqnarray*}
   \pi(u(\lambda_{h_1 + h_2})) &=& \pi(u(\lambda_{h_1})) \, \pi(u(\lambda_{h_2}))
                                =  v(\sigma_{\phi(h_1)}) \, v(\sigma_{\phi(h_2)}) \\
                               &=& v(\sigma_{\phi(h_1)} + \sigma_{\phi(h_2)})
                                =  v(\sigma_{\phi(h_1) + \phi(h_2)}) \quad \mod{\mathbb{T}}.
\end{eqnarray*}
On the other hand, we get $\pi(u(\lambda_{h_1 + h_2})) = v(\sigma_{\phi(h_1 + h_2)}) \mod{\mathbb{T}}$. Since $\{v(\sigma)\}$ are linearly independent, we get $\sigma_{\phi(h_1 + h_2)} = \sigma_{\phi(h_1) + \phi(h_2)}$, and hence
\begin{eqnarray*}
\phi(h_1 + h_2) = \phi(h_1) + \phi(h_2).
\end{eqnarray*}
This means that the map $\phi$ is a group homomorphism. The bijectivity of the $*$-isomorphism $\pi$ leads to that of the group homomorphism $\phi = \phi_\pi$. Since $\{\gamma \cdot \lambda_h \ | \ \gamma \in \G, h \in H_a\} \subset \Lambda(H_a)$ generates $\Lambda(H_a)$, we get $\pi(u(\lambda)) = v(\phi \circ \lambda) \mod{\mathbb{T}}$ for $\lambda \in \Lambda(H_a)$.

We prove that the group isomorphism $\phi = \phi_\pi$ satisfies conditions $(\ref{condition; cohomologous})$ and $(\ref{condition; character})$ in the theorem.
For all $h \in H_a$, there exists $c(h) \in \mathbb{T}$ satisfying
\begin{eqnarray*}
   U_h = \pi \left( \overline{\mu_a(h, -h)} \, u(\lambda_h)\right) = c(h) \, \overline{\mu_b(\phi(h), - \phi(h))} \,  v(\sigma_{\phi(h)}).
\end{eqnarray*}
Since $(e_1, \eta) \in \G$ acts on $\Z^2$ as $(e_1, \eta) \cdot e_1 = 0,$ $(e_1, \eta) \cdot 0 = e_1$, we get
\begin{eqnarray*}
   U_h \, \beta_b(e_1, \eta)(U_h) &=& \pi \left(\overline{\mu_a(h, -h)} \, u(\lambda_h) \, \overline{\mu_a(h, -h)} \, u(-\lambda_h)\right)\\
                                    &=& \overline{\mu_a(h, -h)}^2 \, \widetilde{\mu_a}(\lambda_h, - \lambda_h) = 1.
\end{eqnarray*}
The following equation also holds:
\begin{eqnarray*}
   U_h \, \beta_b(e_1, \eta)(U_h)
&=& c(h) \, \overline{\mu_b(\phi(h), -\phi(h))} \, v(\sigma_{\phi(h)}) \ c(h) \, \overline{\mu_b(\phi(h), -\phi(h))} \, v(\sigma_{-\phi(h)})\\
&=& c(h)^2 \, \overline{\mu_b(\phi(h), -\phi(h))}^2 \, \widetilde{\mu_b}(\sigma_{\phi(h)}, -\sigma_{\phi(h)})
 =  c(h)^2.
\end{eqnarray*}
Thus we have $c(h) \in \{1, -1 \}$ for $h \in H_a$.

Since $\xi \cdot e_1 = e_2$, $\xi \cdot e_2 = 0$ and $\xi \cdot 0 = e_1$, we have
\begin{eqnarray*}
& & U_h \, \beta_b(\xi)(U_h) \, \beta_b(\xi^2)(U_h) \\
&=& \pi \left(\overline{\mu_a(h, -h)} \, u(\lambda_h)\right) \, \pi\left(\chi_a(h) \, \overline{\mu_a(h, -h)} \, u(\xi \cdot \lambda_h)\right) \,
              \pi\left(\overline{\mu_a(h, -h)} \, u(\xi^2 \cdot \lambda_h)\right)\\
&=& \chi_a(h).
\end{eqnarray*}
On the other hand, we have the following:
\begin{eqnarray*}
    U_h \, \beta_b(\xi)(U_h) \, \beta_b(\xi^2)(U_h)
&=& c(h) \, \overline{\mu_b(\phi(h), -\phi(h))} \, v(\sigma_{\phi(h)}) \\
& & \quad    c(h) \, \chi_b(\phi(h)) \, \overline{\mu_b(\phi(h), -\phi(h))} \, v(\xi \cdot \sigma_{\phi(h)}) \\
& & \quad \quad          c(h) \, \overline{\mu_b(\phi(h), -\phi(h))} \, v(\xi^2 \cdot \sigma_{\phi(h)})\\
&=& c(h)^3 \, \chi_b(\phi(h))
 = c(h) \, \chi_b(\phi(h)).
\end{eqnarray*}
It follows that
\begin{eqnarray}\label{equation; uniqueness of c}
c(h) = \overline{\chi_b(\phi(h))} \, \chi_a(h)
\end{eqnarray}
and $\chi_b(\phi(h))^2 = \chi_a(h)^2$, for all $h \in H_a$.

We recall that the algebra $L_{\widetilde{\mu_a}}(\oplus_{\Z^2} H_a)$ is canonically identified wit the infinite tensor product $\bigotimes_{\Z^2} L_{\mu_a} (H_a)$. The unitary $\pi^{-1}(U_h) \in N(H_a, \mu_a) \subset \bigotimes_{\Z^2} L_{\mu_a} (H_a)$ can be written as $1 \otimes u_h^* \otimes u_h$, where $1$ is placed on $-e_1 \in \Z^2$, $u_h^*$ is placed on $0$ and $u_h$ is placed on $e_1$. Since $\eta \in \G$ acts on $\Z^2$ as $\eta \cdot e_1 = - e_1$, $\eta \cdot 0 = 0$, the unitary $\pi^{-1}(\beta_b(\eta)(U_g))$ can be written as $u_g \otimes u_g^* \otimes 1$. We have the following equation:
\begin{eqnarray*}
& & U_g \, \beta_b(\eta)(U_h) \, U_g^* \, \beta_b(\eta)(U_h)^* \\
&=& \pi((1 \otimes u_g^* \otimes u_g) (u_h \otimes u_h^* \otimes 1)
            (1 \otimes u_g^* \otimes u_g)^* (u_h \otimes u_h^* \otimes 1)^*)
 = \mu_a^* \mu_a(g, h).
\end{eqnarray*}
The unitary $U_h$ can be written as $c(h)(1 \otimes v_{\phi(h)}^* \otimes v_{\phi(h)}) \in N(H_b, \mu_b) \subset \bigotimes L_{\mu_b}(H_b)$. Here we write $v_{\phi(h)}$ for the unitary in $L_{\mu_b}(H_b)$ corresponding to $\phi(h)$. The unitary $\beta_b(\eta)(U_g)$ can be written as $c(g)(v_{\phi(g)} \otimes v_{\phi(g)}^* \otimes 1)$. Then we get
\begin{eqnarray*}
& & U_g \, \beta_b(\eta)(U_h) \, U_g^* \, \beta_b(\eta)(U_h^*) \\
&=& (1 \otimes v_{\phi(g)}^* \otimes v_{\phi(g)}) (v_{\phi(h)} \otimes v_{\phi(h)}^* \otimes 1)
            (1 \otimes v_{\phi(g)}^* \otimes v_{\phi(g)})^* (v_{\phi(h)} \otimes v_{\phi(h)}^* \otimes 1)^*\\
&=& \mu_b^* \mu_b(\phi(g), \phi(h)).
\end{eqnarray*}
Thus we get $\mu_a^* \mu_a(g, h) = \mu_b^* \mu_b(\phi(g), \phi(h))$, for all $g, h \in H_a$.
We proved that the group isomorphism $\phi = \phi_\pi$ satisfies conditions $(1)$, $(2)$ and $(3)$.
\end{proof}

From a group homomorphism which satisfies conditions $(2)$ and $(3)$, we construct a $\ast$-homomorphism from $N(H_a, \mu_a)$ to $N(H_b, \mu_b)$ with condition $(\ref{condition; phi circ})$. In the construction, the function $\widehat\mu$ on $\Lambda(H)$ given below is useful. We fix an index for $\Z^2$ as $\Z^2 = \{ k_0, k_1, k_2, \cdots \}$ throughout the rest of this section. For a scalar $2$-cocycle $\mu$ of $H$, we define the function $\widehat\mu$ by
\begin{eqnarray*}
   \widehat\mu(\lambda) &=& \prod_{j = 1}^{n} \mu \left(\sum_{i = 0}^{j-1}\lambda(k_i), \lambda(k_j)\right), \quad \lambda \in \Lambda(H),
\end{eqnarray*}
where $\lambda$ is supported on $\{ k_0, k_1, k_2, \cdots ,k_n \}$. This definition depends on the choice of an order on $\Z^2$. Since $\sum_i \lambda(k_i) = 0$, the function $\widehat\mu$ is also given by the following relation in $\C_\mu(H)$:
\begin{eqnarray*}
   \widehat\mu(\lambda) 1 = u_{\lambda(k_0)} u_{\lambda(k_1)} u_{\lambda(k_2)} \cdots u_{\lambda(k_n)}, \quad \lambda \in \Lambda(H).
\end{eqnarray*}
If $\mu$ is a coboundary, then the definition of $\widehat\mu$ does not depend on the order on $\Z^2$, since $\C_\mu(H)$ is commutative.

\begin{lemma}\label{Lemma for 2-cocycles}
Let $\mu_0$ be another normalized scalar $2$-cocycle for $H$. Let $\widetilde{\mu_0}$ be the scalar $2$-cocycle on $\Lambda(H) \times \Lambda(H)$ given in the same way as equation $(\ref{equation; widetilde})$ in subsection $\ref{subsection; scalar 2-cocycles}$ and let $\widehat{\mu_0}$ be the function on $\Lambda(H)$ constructed from $\mu_0$ in the above manner. If the scalar $2$-cocycles $\mu$ and $\mu_0$ are cohomologous, then for all $\lambda_1, \lambda_2 \in \Lambda(H)$, we have the equation
\begin{eqnarray*}
	\widetilde\mu( \lambda_1, \lambda_2 ) \, \overline{\widehat\mu( \lambda_1 )} \,
         \overline{\widehat\mu( \lambda_2 )} \, \widehat\mu( \lambda_1 + \lambda_2 )
     = \widetilde{\mu_0}( \lambda_1, \lambda_2 ) \, \overline{\widehat{\mu_0}( \lambda_1 )} \,
         \overline{\widehat{\mu_0}( \lambda_2 )} \, \widehat{\mu_0}( \lambda_1 + \lambda_2 ).
\end{eqnarray*}
\end{lemma}
\begin{proof}
We denote by $\{ \nu(g, h) \}$ the scalar $2$-cocycle $\{\overline{\mu_0(g, h)} \mu(g, h) \}$ of $H$. Since $\nu$ is a $2$-coboundary, there exists $\{c(g)\}_{g \in H} \subset \mathbb{T}$ satisfying $\nu(g, h) = b(g) b(h) \overline{b(g + h)}$. Then the map $\widehat\nu$ becomes $\widehat{\nu}(\lambda) = \prod_{i}b(\lambda(k_i))$. Since
\begin{eqnarray*}
    \widehat{\nu}(\lambda_1) \, \widehat{\nu}(\lambda_2) &=& \prod_{i} b(\lambda_1(k_i)) \, b(\lambda_2(k_i)),\\
    \widetilde{\nu}(\lambda_1, \lambda_2) &=& \prod_{i} b(\lambda_1(k_i)) \, b(\lambda_2(k_i)) \, \overline{b(\lambda_1(k_i) + \lambda_2(k_i))},\\
    \widehat{\nu}(\lambda_1 + \lambda_2) &=& \prod_{i} b(\lambda_1(k_i) + \lambda_2(k_i)),
\end{eqnarray*}
we get $\widehat{\nu}(\lambda_1) \, \widehat{\nu}(\lambda_2) = \widetilde{\nu}(\lambda_1, \lambda_2) \, \widehat{\nu}(\lambda_1 + \lambda_2)$. By the definitions of $\widetilde\mu, \widetilde{\mu_0}, \widehat\mu$ and $\widehat{\mu_0}$, the maps $\widehat\nu$ and $\widetilde\nu$ are given by
\begin{eqnarray*}
   \widehat\nu(\lambda)
         = \widehat\mu(\lambda) \, \overline{\widehat{\mu_0}(\lambda)}, \quad
   \widetilde\nu(\lambda_1, \lambda_2)
         = \widetilde\mu(\lambda_1, \lambda_2) \, \overline{\widetilde{\mu_0}(\lambda_1, \lambda_2)},
\end{eqnarray*}
Thus the desired equality immediately follows.
\end{proof}

\begin{proof}[Proof for the second half of Theorem \ref{thm1}]\ \\
Suppose that there exists a group isomorphism $\phi$ satisfying conditions $(\ref{condition; cohomologous})$ and $(\ref{condition; character})$ in the theorem. We prove that there exists a $*$-isomorphism $\pi = \pi_\phi$ from $N(H_a, \mu_a)$ onto $N(H_b, \mu_b)$ preserving the $\G$-actions with condition $(\ref{condition; phi circ})$.

We define a group homomorphism $c_\phi$ from $H_a$ to $\{1, -1\} \subset \mathbb{T}$ by
\begin{eqnarray*}
   c_\phi(h) = \overline{\chi_b(\phi(h))} \chi_a(h), \quad h \in H_a.
\end{eqnarray*}
Let $\widetilde{c_\phi}$ be the group homomorphism from $\Lambda(H_a)$ to $\{1, -1\} \subset \mathbb{T}$ given by
\begin{eqnarray*}
   \widetilde{c_\phi}(\lambda) = \prod_{k \in \Z^2} c_\phi(\lambda(k))^{\gcd(k)}
   = \prod_{k \in \Z^2} \chi_a(\lambda(k))^{\gcd(k)}  \overline{\chi_b(\phi(\lambda(k)))}^{\gcd(k)}
   , \quad \lambda \in \Lambda(H_a).
\end{eqnarray*}
We define a linear map $\pi$ from the group algebra $\C_{\widetilde\mu_a}(\Lambda(H_a))$ onto $\C_{\widetilde\mu_b}(\Lambda(H_b))$ by
\begin{eqnarray*}
   \pi\left(\overline{\widehat{\mu_a}(\lambda)} \, u(\lambda)\right)
 = \widetilde{c_\phi}(\lambda) \, \overline{\widehat{\mu_b}(\phi \circ \lambda)} \, v(\phi \circ \lambda),
 \quad \lambda \in \Lambda(H_a).
\end{eqnarray*}
By direct computations, for all $\lambda_1, \lambda_2 \in \Lambda(H_a)$, we get
\begin{eqnarray*}
& & \pi\left(\overline{\widehat{\mu_a}(\lambda_1)} \, u(\lambda_1)\right) \, \pi\left(\overline{\widehat{\mu_a}(\lambda_2)} \, u(\lambda_2)\right) \\
&=& \widetilde{c_\phi}(\lambda_1) \, \widetilde{c_\phi}(\lambda_2) \,
    \overline{\widehat{\mu_b}(\phi \circ \lambda_1)} \, \overline{\widehat{\mu_b}(\phi \circ \lambda_2)} \,
    v(\phi \circ \lambda_1) \, v(\phi \circ \lambda_2) \\
&=& \widetilde{c_\phi}(\lambda_1 + \lambda_2) \,
    \overline{\widehat{\mu_b}(\phi \circ \lambda_1)} \, \overline{\widehat{\mu_b}(\phi \circ \lambda_2)} \,
    \widetilde{\mu_b}(\phi \circ \lambda_1, \phi \circ \lambda_2) \, v(\phi \circ (\lambda_1 + \lambda_2)).
\end{eqnarray*}
On the other hand, we have the following equation:
\begin{eqnarray*}
& & \pi\left(\overline{\widehat{\mu_a}(\lambda_1)} \, \overline{\widehat{\mu_a}(\lambda_2)} \, u(\lambda_1) \, u(\lambda_2)\right) \\
&=& \pi\left(\overline{\widehat{\mu_a}(\lambda_1)} \, \overline{\widehat{\mu_a}(\lambda_2)} \,
        \widetilde{\mu_a}(\lambda_1, \lambda_2) \, u(\lambda_1 + \lambda_2)\right) \\
&=& \overline{\widehat{\mu_a}(\lambda_1)} \, \overline{\widehat{\mu_a}(\lambda_2)} \, \widetilde{\mu_a}(\lambda_1, \lambda_2) \,
      \widehat{\mu_a}(\lambda_1 + \lambda_2) \, \pi\left(\overline{\widehat{\mu_a}(\lambda_1 + \lambda_2)} \, u(\lambda_1 + \lambda_2)\right) \\
&=& \widetilde{c_\phi}(\lambda_1 + \lambda_2) \,
    \overline{\widehat{\mu_a}(\lambda_1)} \, \overline{\widehat{\mu_a}(\lambda_2)} \, \widetilde{\mu_a}(\lambda_1, \lambda_2)          \,  \widehat{\mu_a}(\lambda_1 + \lambda_2)\\
& & \quad \overline{\widehat{\mu_b}(\phi \circ(\lambda_1 + \lambda_2))} \, v(\phi \circ (\lambda_1 + \lambda_2)).
\end{eqnarray*}
By Lemma \ref{Lemma for 2-cocycles} and condition $(\ref{condition; cohomologous})$ for the group isomorphism $\phi$ in the theorem, we have that
\begin{eqnarray*}
& & \overline{\widehat{\mu_b}(\phi \circ \lambda_1)} \, \overline{\widehat{\mu_b}(\phi \circ \lambda_2)} \,
    \widetilde{\mu_b}(\phi \circ \lambda_1, \phi \circ \lambda_2)\\
&=& \overline{\widehat{\mu_a}(\lambda_1)} \, \overline{\widehat{\mu_a}(\lambda_2)} \, \widetilde{\mu_a}(\lambda_1, \lambda_2)          \, \widehat{\mu_a}(\lambda_1 + \lambda_2) \,
    \overline{\widehat{\mu_b}(\phi \circ (\lambda_1) + \phi \circ (\lambda_2))}.
\end{eqnarray*}
Therefore we get $\pi(u(\lambda_1)) \, \pi(u(\lambda_2)) = \pi(u(\lambda_1) \, u(\lambda_2))$. The linear map $\pi$ also preserves the $*$-operation. As a consequence, $\pi$ is a $*$-isomorphism from $\C_{\widetilde\mu_a}(\Lambda(H_a))$ onto $\C_{\widetilde\mu_b}(\Lambda(H_b))$ and this preserves the trace. The map $\pi = \pi_\phi$ is extended to a normal $*$-isomorphism from $N(H_a, \mu_a)$ onto $N(H_b, \mu_b)$.

We next prove that this $\pi$ preserves the $\G$-actions. The group homomorphism $\widetilde{c_\phi}$ from $\Lambda(H_a)$ to $\{ 1, -1\}$ is invariant under the action of ${\rm SL}(2, \Z)$, by Lemma \ref{det and gcd} $(1)$. The scalar $2$-cocycle $\nu(g, h) = \mu_a(g, h) \overline{\mu_b(\phi(h), \phi(g))}$
satisfies $\nu(g, h) = \nu(h, g)$ by condition $(\ref{condition; cohomologous})$, so the function $\widehat\nu(\cdot) = \widehat{\mu_a}(\cdot) \overline{\widehat{\mu_b}(\phi \circ \cdot)}$ on $\Lambda(H_a)$ does not depend on the order on $\Z^2$ chosen before. Since
\begin{eqnarray*}
    \pi \circ \beta_a(\gamma) (u(\lambda))
&=& \pi(u(\gamma \cdot \lambda))\\
&=& \widetilde{c_\phi}(\gamma \cdot \lambda)
    \widehat{\mu_a}(\gamma \cdot \lambda) \overline{\widehat{\mu_b}(\phi \circ (\gamma \cdot \lambda))}
    v(\phi \circ (\gamma \cdot \lambda)) \\
&=& \widetilde{c_\phi}(\lambda)
    \widehat{\mu_a}(\lambda) \overline{\widehat{\mu_b}(\phi \circ \lambda)}
    v(\gamma \cdot (\phi \circ \lambda))\\
&=& \beta_b(\gamma)\left(\widetilde{c_\phi}(\lambda)
    \widehat{\mu_a}(\lambda) \overline{\widehat{\mu_b}(\phi \circ \lambda)}
    v(\phi \circ \lambda)\right)\\
&=& \beta_b(\gamma) \circ \pi(u(\lambda)), \quad \gamma \in {\rm SL}(2, \Z), \lambda \in \Lambda(H_a),
\end{eqnarray*}
it turns out that the $*$-isomorphism $\pi$ preserves the ${\rm SL}(2, \Z)$-action.

For all $\lambda \in \Lambda(H_a)$ and $k \in \Z^2$, we have
\begin{eqnarray*}
& & \pi \circ \beta_a(k) (u(\lambda))
 =  \pi \left( \prod_{l \in \Z^2} \chi_a(\lambda(l))^{\det(k, l)} u(k \cdot \lambda) \right)\\
&=& \prod_{l \in \Z^2} c_{\phi}((k \cdot \lambda)(l))^{\gcd(l)}
    \prod_{l \in \Z^2} \chi_a(\lambda(l))^{\det(k, l)}
    \widehat{\mu_a}(k \cdot \lambda) \, \overline{\widehat{\mu_b}(\phi \circ (k \cdot \lambda))} \,
        v(\phi \circ (k \cdot \lambda)).
\end{eqnarray*}
Since
$c_\phi(h)^{\det(k, l)} c_\phi(h)^{\gcd(k + l)} = c_\phi(h)^{\gcd(k)} c_\phi(h)^{\gcd(l)}$,
by Lemma \ref{det and gcd} $(2)$, the unitary $\pi \circ \beta_a(k) (u(\lambda))$ equals to
\begin{eqnarray*}
& & \prod_{l \in \Z^2} c_{\phi}(\lambda(l))^{\gcd(k + l)}
    \prod_{l \in \Z^2} \left(c_{\phi}(\lambda(l))^{\det(k, l)} \chi_b(\phi \circ \lambda(l))^{\det(k, l)}\right)
    \\
& & \qquad \widehat{\mu_a}(k \cdot \lambda) \, \overline{\widehat{\mu_b}(\phi \circ (k \cdot \lambda))} \,
        v(\phi \circ (k \cdot \lambda))\\
&=& \prod_{l \in \Z^2} c_{\phi}(\lambda(l))^{\gcd(k)}
    \prod_{l \in \Z^2} c_{\phi}(\lambda(l))^{\gcd(l)} \prod_{l \in \Z^2} \chi_b(\phi \circ \lambda(l))^{\det(k, l)}
    \\
& & \qquad \widehat{\mu_a}(k \cdot \lambda) \, \overline{\widehat{\mu_b}(\phi \circ (k \cdot \lambda))} \,
        v(k \cdot (\phi \circ \lambda))\\
&=& \widetilde{c_{\phi}}(\lambda) \,
    \widehat{\mu_a}(k \cdot \lambda) \, \overline{\widehat{\mu_b}(k \cdot (\phi \circ \lambda))} \,
    \beta_b(k) (v(\phi \cdot \lambda))\\
&=& \beta_b(k) \circ \pi (u(\lambda)).
\end{eqnarray*}
This means that the $*$-isomorphism $\pi$ preserves the $\Z^2$-actions.

We get the $*$-isomorphism $\pi = \pi_\phi$ from $N(H_a, \mu_a)$ onto $N(H_b, \mu_b)$ giving conjugacy between $\beta_a$ and $\beta_b$.
\end{proof}

\begin{remark}\label{trace preserving}
The proof of the first half of Theorem \ref{thm1} shows that any isomorphism $\pi$ giving conjugacy between $\beta_a$ and $\beta_b$ is of the form $\pi_\phi$. This means that an isomorphism which gives conjugacy between two twisted Bernoulli shift actions must be trace preserving.
\end{remark}

This proof shows that an isomorphism giving conjugacy between the two actions $\beta(H_a, \mu_a, \chi_a)$, $\beta(H_b, \mu_b, \chi_b)$ is of a very special form derived from a group isomorphism between $H_a$ and $H_b$. Taking notice of this fact, we can describe the centralizer of a twisted Bernoulli shift action. We define two topological groups before we state Theorem \ref{thm2}.

Let $\beta$ be a trace preserving action of some group $\Gamma$ on a separable finite von Neumann algebra $(N, \tr)$. We denote by ${\rm Aut}(N, \beta)$ the group of all automorphisms which commute with the action $\beta$, that is,
\begin{eqnarray*}
\{ \alpha \in {\rm Aut}(N) \ |
             \ \beta(\gamma) \circ \alpha =  \alpha \circ \beta(\gamma), \quad \gamma \in \Gamma \}.
\end{eqnarray*}
We regard the group ${\rm Aut}(N, \beta)$ as a topological group equipped with the pointwise-strong topology. When $\beta$ is a twisted Bernoulli shift action on $N$, an automorphism $\alpha$ commuting with $\beta$ is necessarily trace preserving by Remark \ref{trace preserving}. We consider that ${\rm Aut}(N, \beta)$ is equipped with the pointwise-$2$-norm topology.

Let ${\rm Aut}(H, \mu, \chi)$ be the group of all automorphisms of an abelian group $H$ which preserve its $2$-cocycle $\mu$ and character $\chi$, that is,
\begin{eqnarray*}
   \{ \phi \in {\rm Aut}(H) \ |
             \ \mu(g, h) = \mu(\phi(g), \phi(h)), \chi(g) = \chi(\phi(g)), \quad g, h \in H \}.
\end{eqnarray*}
We define the topology of ${\rm Aut}(H, \mu, \chi)$ by pointwise convergence.

\begin{theorem}\label{thm2}
For $\pi \in {\rm Aut}(N(H, \mu), \beta(H, \mu, \chi))$, there exists a unique element $\phi = \phi_\pi \in {\rm Aut}(H, \mu^*\mu, \chi^2)$ satisfying $\pi(u(\lambda)) = u(\phi \circ \lambda)\ \mathrm{mod} \ \mathbb{T}$ for $\lambda \in \Lambda(H)$. The map $\pi \mapsto \phi_\pi$ gives an isomorphism between two topological groups
\begin{eqnarray*}
   {\rm Aut}(N(H, \mu), \beta(H, \mu, \chi)) \cong {\rm Aut}(H, \mu^*\mu, \chi^2).
\end{eqnarray*}
\end{theorem}

\begin{proof}
We use the notations in the proof of the previous theorem letting $H_a = H_b = H$, $\mu_a = \mu_b = \mu$ and $\chi_a = \chi_b = \chi$. Denote $N = N(H, \mu)$ and $\beta = \beta(H, \mu, \chi)$. We have already have shown the first claim. Let
\begin{eqnarray*}
{\rm Aut}(H, \mu^*\mu, \chi^2) \ni \phi &\mapsto& \pi_\phi \in {\rm Aut}(N, \beta),
\end{eqnarray*}
be the map given as in the proof of Theorem \ref{thm1}, that is,
\begin{eqnarray*}
   \pi_\phi \left(\overline{\widehat{\mu}(\lambda)} \, u(\lambda)\right)
 = \widetilde{c_\phi}(\lambda) \, \overline{\widehat{\mu}(\phi \circ \lambda)} \, u(\phi \circ \lambda), \quad \lambda \in \Lambda(H),
\end{eqnarray*}
where
\begin{eqnarray*}
   \widetilde{c_\phi}(\lambda) = \prod_{k \in \Z^2} \chi(\lambda(k))^{\gcd(k)} \, \overline{\chi(\phi \circ \lambda(k))}^{\gcd(k)}.
\end{eqnarray*}

It is easy to prove that $\phi_{\pi_\phi} = \phi$ by the definition. Thus the map $\pi \mapsto \phi_\pi$ is surjective.
This map is also injective. Let $\phi$ be an element of ${\rm Aut}(H, \mu^*\mu, \chi^2)$. Suppose that $\pi$ is an arbitrary element of ${\rm Aut}(N,\beta)$ satisfying $\phi = \phi_\pi$. The set $\{ \beta(\gamma) (u(\lambda_h)) \ | \ h \in H, \ \gamma \in \G\}$ generates $N$, so we have only to prove the uniqueness of $c(h) \in \mathbb{T}$ satisfying
\begin{eqnarray*}
    \pi \left(\overline{\mu(h, -h)} \, u(\lambda_h)\right) = c(h) \, \overline{\mu(\phi(h), -\phi(h))} \, u(\lambda_{\phi(h)}),
\end{eqnarray*}
for all $h \in H$. In the proof of the first half of the previous theorem (equation $(\ref{equation; uniqueness of c})$), we have already shown that $c(h) = \chi(h) \overline{\chi(\phi(h))}$. Thus the $*$-isomorphism $\pi$ is uniquely determined and the map $\pi \mapsto \phi_\pi$ is injective.

We prove the two maps $\phi \mapsto \pi_\phi$ and $\pi \mapsto \phi_\pi$ are continuous. Let $(\phi_i)$ be a net in ${\rm Aut}(H, \mu^*\mu, \chi^2)$ converging to $\phi$. For all $h \in H$, we have
\begin{eqnarray*}
    \pi_{\phi_i} \left(\overline{\chi(h)} \, \overline{\mu(h, -h)} \, u(\lambda_h)\right) = \overline{\chi(\phi_i(h))} \, \overline{\mu(\phi_i(h), -\phi_i(h))} \, u(\lambda_{\phi_i(h)}).
\end{eqnarray*}
The right side of the equation converges to
\begin{eqnarray*}
    \overline{\chi(\phi(h))} \, \overline{\mu(\phi(h), -\phi(h))} \, u(\lambda_{\phi(h)}) = \pi_\phi \left(\overline{\chi(h)} \, \overline{\mu(h, -h)} \, u(\lambda_h)\right).
\end{eqnarray*}
This proves that $\pi_{\phi_i}$ converges to $\pi_\phi$ in pointwise 2-norm topology on the generating set $\{ \beta(\gamma) (u(\lambda_h)) \ | \ h \in H, \ \gamma \in \G\}$ of $N$. Thus $\pi_{\phi_i}$ converges to $\pi_\phi$ on $N$.

Conversely, let $(\pi_i)$ be a net in ${\rm Aut}(N, \beta)$ converging to $\pi$. For all $h \in H$, we get $\pi_i (u(\lambda_h)) = u(\lambda_{\phi_{\pi_i}(h)}) \mod{\mathbb{T}}$. The left side of the equation converges to $\pi (u(\lambda_h)) = u(\lambda_{\phi_\pi(h)})$. If $\phi_{\pi_i}(h) \neq \phi_{\pi}(h)$, then the distance between $\mathbb{T} u(\lambda_{\phi_\pi(h)}) $ and $\mathbb{T} u(\lambda_{\phi_{\pi_i}(h)})$ is $\sqrt{2}$ in the $2$-norm. Thus $\phi_{\pi_i}(h) = \phi_\pi(h)$ for large enough $i$. This means that $(\phi_{\pi_i})$ converges to $\phi_\pi$.

As a consequence, the two maps $\phi \mapsto \pi_\phi$ and $\pi \mapsto \phi_\pi$ are continuous group homomorphisms and inverse maps of each other.
\end{proof}

\section{Examples}\label{examples}
\subsection{Twisted Bernoulli shift actions on $L^\infty(X)$}
In this subsection, we consider the case of $\mu = 1$ and $H \neq \{ 0 \}$. Then the algebra $N(H, 1)$ is abelian and has a faithful normal state, so it is isomorphic to $L^\infty(X)$, where $X$ is a standard probability space. The measure of $X$ is determined by the trace on $N(H, 1)$. Furthermore, $X$ is non-atomic, since $N(H, 1)$ is infinite dimensional and the action $\beta(H, 1, \chi)$ is ergodic. As corollaries of Theorems \ref{thm1} and \ref{thm2}, we get trace preserving $\G$-actions on $L^\infty(X)$ whose centralizers are isomorphic to some prescribed groups.

\begin{remark}
The $\G$-action on $X$ defined by $\beta = \beta(H, 1, \chi)$ is free. An automorphism $\alpha \in {\rm Aut}(L^\infty(X), \beta)$ is free or the identity map for any twisted Bernoulli shift action $\beta$ on $L^\infty(X)$. This is proved as follows. We identify $\beta(\gamma)$ $(\gamma \in \G)$ and $\alpha$ with measure preserving Borel isomorphisms on $X$ here. Suppose that there exists a non-null Borel subset $Y \subset X$ whose elements are fixed under $\alpha$. All elements in $\widetilde{Y} = \cup\{ \beta(\gamma)(Y) \ | \ \gamma \in \G \}$ are fixed under $\alpha$. By the ergodicity of $\beta$, the measure of $\widetilde{Y}$ is $1$. Then $\alpha$ is the identity map of $L^\infty(X)$.
\end{remark}

\begin{corollary}\label{Aut(H)}
For any abelian countable discrete group $H \neq \{ 0 \}$, there exists a trace preserving essentially free ergodic action $\beta$ of $\Z^2 \rtimes {\rm SL}(2, \Z)$ on $L^\infty(X)$ satisfying
\begin{eqnarray*}
   {\rm Aut}(L^\infty(X), \beta) \cong {\rm Aut}(H).
\end{eqnarray*}
\end{corollary}

\begin{proof}
When we define $\beta = \beta(H, 1, 1)$, we have the above relation by Theorem \ref{thm2}.
\end{proof}

In the next corollary we use the effect of twisting by a character $\chi$.

\begin{corollary}
For every abelian countable discrete group $H \neq \{ 0 \}$, there exist continuously many trace preserving essentially free ergodic actions $\{ \beta_c \}$ of $\Z^2 \rtimes {\rm SL}(2, \Z)$ on $L^\infty(X)$ which are mutually non-conjugate and satisfy
\begin{eqnarray*}
   {\rm Aut}(L^\infty(X), \beta_c) \cong H \rtimes {\rm Aut}(H).
\end{eqnarray*}
Here the topology of $H \rtimes {\rm Aut}(H)$ is the product of the discrete topology on $H$ and the pointwise convergence topology on ${\rm Aut}(H)$.
\end{corollary}

\begin{proof}
Let  $c \in \{ e^{i \pi t} \ | \ t \in (0, 1/2) \setminus \Q\}$. We put $\beta_c = \beta(H \oplus \Z, 1, 1 \times \chi_c)$, where the character $\chi_c$ of $\Z$ is defined as $\chi_c(n) = c^n$. By Theorem \ref{thm2}, we get
\begin{eqnarray*}
   {\rm Aut}(L^\infty(X), \beta_c) \cong {\rm Aut}(H \oplus \Z, 1, 1 \times \chi_{c}^2).
\end{eqnarray*}
Since the character $\chi_{c}^2$ is injective, a group automorphism $\alpha \in {\rm Aut}(H \oplus \Z, 1, 1 \times \chi_{c}^2)$ preserves the second entry. For all $\alpha \in {\rm Aut}(H \oplus \Z, 1, 1 \times \chi_{c}^2)$, there exist $\phi_\alpha \in {\rm Aut}(H)$ and $h_\alpha \in H$ satisfying
\begin{eqnarray*}
   \alpha(h, n) = (\phi_\alpha(h) + n h_\alpha, n), \quad (h, n) \in H \oplus \Z.
\end{eqnarray*}
The map ${\rm Aut}(H \oplus \Z, 1, 1 \times \chi_{c}^2) \ni \alpha \mapsto (h_\alpha, \phi_\alpha) \in H \rtimes {\rm Aut}(H)$ is a homeomorphic group isomorphism.

If $c_1, c_2 \in \{ e^{i \pi t} \ | \ t \in (0, 1/2) \setminus \Q\}$ and $c_1 \neq c_2$, then there exists no isomorphism from $H \oplus \Z$ to $H \oplus \Z$ whose pull back of the character $1 \times \chi_{c_2}^2$ is equal to $1 \times \chi_{c_1}^2$. The two actions $\beta_{c_1}$ and $\beta_{c_2}$ are not conjugate by Theorem \ref{thm1}.
\end{proof}

\begin{corollary}\label{Corollary; trivial centralizer}
There exist continuously many trace preserving essentially free ergodic actions $\{\beta_c\}$ of $\Z^2 \rtimes {\rm SL}(2, \Z)$ on $L^\infty(X)$ which are mutually non-conjugate and have the trivial centralizer
${\rm Aut}(L^\infty(X), \beta_c) = \{ \id_{L^\infty X} \}$.
\end{corollary}

\begin{proof}
Let $\{ \chi_c \ | \ c = e^{i \pi t}, t \in (0, 1/2) \}$ be characters of $\Z$ such that $\chi_c(m) = c^m$. Since $\chi_c^2(1)$ is in the upper half plane, the identity map is the only automorphism of $\Z$ preserving $\chi_c^2$. By Theorem \ref{thm2}, we get ${\rm Aut}(\beta(\Z, 1, \chi_c)) = \{ \id \}$.

If $\beta(\Z, 1, \chi_{c_1}), \beta(\Z, 1, \chi_{c_2})$ are conjugate, then there exists a group isomorphism on $\Z$ whose pull back of $\chi_{c_2}^2$ is $\chi_{c_1}^2$ by Theorem \ref{thm1}. This means $c_1 = c_2$. Thus the actions $\{ \beta(\Z, 1, \chi_c) \}$ are mutually non-conjugate.
\end{proof}

\subsection{Twisted Bernoulli shift actions on the AFD factor of type ${\rm II}_1$}
Firstly, we find a condition that the finite von Neumann algebra $N(H, \mu)$ is the AFD factor of type ${\rm II}_1$.
\begin{lemma}\label{AFD II_1}
For an abelian countable discrete group $H \neq \{ 0 \}$ and its normalized scalar $2$-cocycle $\mu$, the following statements are equivalent:
\begin{enumerate}
\item
The algebra $N(H, \mu)$ is the AFD factor of type ${\rm II}_1$.
\item
The group von Neumann algebra $L_\mu(H)$ twisted by the scalar $2$-cocycle $\mu$ is a factor (of type $\mathrm{II}_1$ or $\mathrm{I}_n$).
\item
For all $g \in H \setminus \{ 0 \}$, there exists $h \in H$ such that $\mu(g, h) \neq \mu(h, g)$.
\end{enumerate}
\end{lemma}

\begin{proof}
The amenability of the group $\Lambda(H)$ leads the injectivity for $N(H, \mu)$. The injectivity for $N(H, \mu)$ implies that $N(H, \mu)$ is approximately finite dimensional (\cite{Co}). We have only to show the equivalence of conditions $(2)$, $(3)$ and
\begin{center}
$(1)'$ \ The algebra $N(H, \mu)$ is a factor.
\end{center}
By using Fourier expansion it is easy to see that condition $(2)$ holds true if and only if for any $g \in H \setminus \{0\}$ there exists $h \in H$ satisfying $ u_g u_h \neq u_h u_g$. This is equivalent to condition $(3)$. Similarly, condition $(1)'$ is equivalent to
\begin{center}
$(1)''$ \ For any $\lambda_1 \in \Lambda(H) \setminus \{0\}$, there exists $\lambda_2 \in \Lambda(H)$ satisfying
\begin{eqnarray*}
\overline{\widetilde\mu(\lambda_2, \lambda_1)} \widetilde\mu(\lambda_1, \lambda_2) \neq 1.
\end{eqnarray*}
\end{center}
Suppose condition $(3)$. For any $\lambda_1$, choose element $k, l \in \Z^2$ so that $k \in \textrm{supp}(\lambda_1)$ and $l \notin \textrm{supp}(\lambda_1)$. By condition $(3)$, there exists $h \in H$ satisfying $\mu^* \mu(\lambda_1(k), h) \neq 1$. Let $\lambda_2$ be the element in $\Lambda(H)$ which takes $h$ at $k$, $-h$ at $l$ and $0$ for the other places. The element $\lambda_2$ satisfies $\overline{\widetilde\mu(\lambda_2, \lambda_1)} \widetilde\mu(\lambda_1, \lambda_2) = \mu^* \mu(\lambda_1(k), h) \neq 1$. Here we get condition $(1)''$. The implication from $(1)''$ to $(3)$ is easily shown.
\end{proof}
\begin{remark}
The twisted Bernoulli shift action $\beta = \beta(H, \mu, \chi)$ is an outer action of $\G$. Any non-trivial automorphism in ${\rm Aut}(R, \beta(H, \mu, \chi))$ is also outer. This is proved by the weak mixing property of the action $\beta(H, \mu, \chi)$ as follows. If $\alpha \in {\rm Aut}(R, \beta(H, \mu, \chi))$ is an inner automorphism ${\rm Ad}(u)$, then we have
\begin{eqnarray*}
     {\rm Ad}(\beta(\gamma)(u))(x)
 &=& \beta(\gamma) (u \beta(\gamma)^{-1}(x) u^*) = \beta(\gamma) \circ \alpha \circ \beta(\gamma)^{-1}(x) \\
 &=& \alpha (x) = {\rm Ad}(u) (x),
\end{eqnarray*}
for all $x \in R$ and $\gamma \in \G$. Since $\mathrm{Ad}(\beta(\gamma)(u)u^*) = \mathrm{id}$, $\C u \subset R$ is an invariant subspace of the action $\beta$. The only subspace invariant under the weakly mixing action $\beta$ is $\C 1$ (Proposition \ref{Weakly Mixing}), thus we get $\alpha = \id$.
\end{remark}

Using Theorems \ref{thm1} and \ref{thm2}, we give continuously many actions of $\G$ on $R$ such that there exists no commuting automorphism except for trivial one.

\begin{corollary}
There exist continuously many ergodic outer actions $\{ \beta_c \}$ of $\Z^2 \rtimes {\rm SL}(2, \Z)$ on the AFD factor $R$ of type $\mathrm{II}_1$ which are mutually non-conjugate and have the trivial centralizer $\mathrm{Aut}(R, \beta_c) = \{\mathrm{id}_R\}$.
\end{corollary}

\begin{proof}
We can choose and fix a character $\chi$ on $\Z^2$ such that $\chi^2$ is injective. Let $\{ \mu_c \ | \ c = e^{i \pi t}, t \in (0, 1/2) \setminus \Q\}$ be scalar $2$-cocycles for $\Z^2$ defined by
\begin{eqnarray*}
   \mu_c
   \left(
   \left(
   \begin{array}{c}
       s_1 \\
       t_1
   \end{array}
   \right),
   \left(
   \begin{array}{c}
       s_2 \\
       t_2
   \end{array}
   \right)
   \right)
   = c^{s_1 t_2 - t_1 s_2}, \quad s_1, t_1, s_2, t_2 \in \Z.
\end{eqnarray*}
We put $\beta_c = \beta(\Z^2, \mu_c, \chi)$. The $2$-cocycle $\mu_c$ satisfies condition $(3)$ in Lemma \ref{AFD II_1}. Thus $\beta_c$ defines a $\G$-action on $R$. By Theorem \ref{thm2}, we get the following isomorphism between topological groups:
\begin{eqnarray*}
   {\rm Aut}(R, \beta_c) \cong {\rm Aut}(\Z^2, \mu_c^* \mu_c, \chi^2) = {\rm Aut}(\Z^2, \mu_{c^2}, \chi^2).
\end{eqnarray*}
Since the character $\chi^2$ of $\Z^2$ is injective, so the group of the right side is $\{ \id |_{\Z^2} \}$. This means that the action $\beta_c$ has trivial centralizers.

Finally, we prove that the actions $\{ \beta_c \ | \ c = e^{i \pi t}, t \in (0, 1/2) \setminus \Q\}$ are mutually non-conjugate. Suppose that actions $\beta_{c_1}$ and $\beta_{c_2}$ are conjugate. By Theorem \ref{thm1}, there exists a group isomorphism $\phi$ of $\Z^2$ satisfying
\begin{eqnarray*}
   \mu_{c_1^2}(g, h)
 = \mu_{c_2^2}(\phi(g), \phi(h)), \quad g, h \in \Z^2.
\end{eqnarray*}
A group isomorphism of $\Z^2$ is given by an element of $\mathrm{GL}(2, \Z)$. If the automorphism $\phi$ is given by an element of ${\rm SL}(2, \Z)$, we get $c_2^2 = c_1^2$. If $\phi$ is given by an element of $\mathrm{GL}(2, \Z) \setminus {\rm SL}(2, \Z)$, then we get $c_2^2 = - c_1^2$. Since both $c_1$ and $c_2$ have the form $e^{i \pi t}, t \in (0, 1/2)$, we get $c_1 = c_2$.
\end{proof}

Any cyclic group of an odd order can be realized as the centralizer of a twisted Bernoulli shift actions on $R$.

\begin{corollary}\label{Odd Prime numbers}
Let $q$ be an odd natural number $\ge 3$ and denote by $H_q$ the abelian group $(\Z/q\Z)^2$. We define the $2$-cocycle $\mu_q$ and the character $\chi_q$ on $H_q$ as
\begin{eqnarray*}
\begin{array}{cc}
\mu_q \left(
    \left(
         \begin{array}{c}
            s_1 \\
            t_1
         \end{array}
         \right),
    \left(
         \begin{array}{c}
            s_2 \\
            t_2
         \end{array}
    \right)\right)
    =
    \exp{( 2 \pi i s_1 t_2/ q)}, &
\chi_q
    \left(\left(
         \begin{array}{c}
            s_1 \\
            t_1
         \end{array}
    \right)\right)
    =
    \exp{(2 \pi i s_1 / q)}.
\end{array}
\end{eqnarray*}
Then the algebra $N(H_q, \mu_q)$ is the AFD factor $R$ of type $\mathrm{II}_1$ and
the centralizer of the twisted Bernoulli shift action $\beta_a = \beta(H_q, \mu_q, \chi_q)$ is isomorphic to $\Z / q \Z$.
\end{corollary}

\begin{proof}
By Lemma \ref{AFD II_1}, it is shown that the algebra $N(H_Q, \mu_Q)$ is the AFD factor of type ${\rm II}_1$. Using Theorem \ref{thm2}, we have only to prove that ${\rm Aut}(H_q, \mu_q^* \mu_q, \chi_q^2) \cong \Z/q\Z$.

Let $\phi$ be in ${\rm Aut}(H_q, \mu_q^* \mu_q, \chi_q^2)$. The automorphism $\phi$ of $H_q$ is given by a $2 \times 2$ matrix $A$ of $\Z / q\Z$. Since $\phi$ preserves $\mu_q^* \mu_q$, the determinant of $A$ must be $1$. Since $q$ is odd, the value of $\chi_q^2$ determines the first entry of $(\Z / q\Z)^2$ and $\phi$ preserves $\chi_q^2$. The matrix $A$ is of the form
\begin{eqnarray*}
   \left(
   \begin{array}{cc}
        1 & 0 \\
        t_{\phi} & 1
   \end{array}
   \right), \quad t_\phi \in \Z / q\Z.
\end{eqnarray*}
The map $\phi \mapsto t_\phi$ is an isomorphism. In turn, if the matrix $A$ is of this form, it defines an element in ${\rm Aut}(H_q, \mu_q^* \mu_q, \chi_q^2)$.
\end{proof}

\begin{corollary}\label{Collection of Odd Prime numbers}
For a set $Q$ consisting of odd prime numbers, let $\beta_Q$ be the tensor product $\bigotimes_{q \in Q} \beta_q$ of the actions $\beta_q$ on the AFD factor of type $II_1$. The centralizer of $\beta_Q$ is isomorphic to $\prod_{q \in Q} \Z/q\Z$.
\end{corollary}

\begin{proof}
The action $\beta_Q$ is the twisted Bernoulli shift action $\beta(H_Q, \mu_Q, \chi_Q)$, where $H_Q$ is the abelian group $\oplus_{q \in Q}H_q$ and the scalar $2$-cocycle $\mu_q$ and a character $\chi_Q$ on $H_Q$ are given by
\begin{eqnarray*}
   \mu_Q((s_q), (t_q)) &=& \prod_{q \in Q} \mu_q (s_q, t_q),\\
   \chi_Q((s_q)) &=& \prod_{q \in Q} \chi_q (s_q), \quad (s_q), (t_q) \in H_Q, s_q, t_q \in H_q.
\end{eqnarray*}
Using Theorem \ref{thm2}, we have only to prove
\begin{eqnarray*}
   {\rm Aut}(H_Q, \mu_Q^* \mu_Q, \chi_Q^2) \cong
           \prod_{q \in Q} \Z/q\Z.
\end{eqnarray*}
A group automorphism $\phi$ of $H_Q = \oplus_{q \in Q} H_q$ has a form $\phi( (k_q) ) = ( \phi_q(k_q) )$, for some $\{ \phi_q \in {\rm Aut}(H_q)\}$. Thus we get
\begin{eqnarray*}
   {\rm Aut}(H_Q, \mu_Q^* \mu_Q, \chi_Q^2) \cong
   \prod_{q \in Q} {\rm Aut}(H_q, \mu_q^* \mu_q, \chi_q^2).
\end{eqnarray*}
Together with the previous corollary, we get the conclusion.
\end{proof}

\begin{remark}
If $Q_1 \neq Q_2$, then the two groups $\prod_{q \in Q_1} \Z / q\Z$ and $\prod_{q \in Q_2} \Z / q\Z$ are not isomorphic. The continuously many outer actions $\{ \beta_Q\}$ are distinguished in view of conjugacy only by using the centralizers $\{ {\rm Aut}(R, \beta_Q) \}$.
\end{remark}

\section{Malleability and rigidity arguments}\label{Section; outer conjugacy}
In this section, we give {\it malleability} and {\it rigidity} type arguments invented by S. Popa, in order to examine weak $1$-cocycles for actions. See Popa \cite{Po2}, \cite{Po3}, \cite{Po4} and Popa--Sasyk \cite{PoSa} for the references. S. Popa in \cite{Po3} showed that every $1$-cocycle for a Connes-St{\o}rmer Bernoulli shift by property (T) group (or $w$-rigid group like $\G$) vanishes modulo scalars. As a consequence, two such actions are cocycle conjugate if and only if they are conjugate. In our case, $1$-cocycles do not vanish modulo scalars but they are still in the situation that cocycle (outer) conjugacy implies conjugacy. We need the following notion to examine outer conjugacy of two group actions.

\begin{definition}
Let $\alpha$ be an action of discrete group $\Gamma$ on a von Neumann algebra $\mathcal{M}$. A {\rm weak\ $1$-cocycle} for $\alpha$ is a map $w : \Gamma \rightarrow \mathcal{U}(M)$ satisfying
\begin{eqnarray*}
    w_{g h} = w_g \alpha_g (w_h) \quad {\rm mod} \ \mathbb{T}, \quad g, h \in \Gamma.
\end{eqnarray*}
The weak $1$-cocycle $w$ is called a {\rm weak $1$-coboundary} if there exists a unitary $v \in \mathcal{U}(M)$ satisfying $w_g = v \alpha_g(v)^* \ {\rm mod}\ \mathbb{T}$. Two weak $1$-cocycles $w$ and $w^\prime$ are said to be {\rm equivalent} when $w^\prime_g = v w_g \alpha_g(v)^*  \ {\rm mod}\ \mathbb{T}$ for some $v \in \mathcal{U}(M)$.
\end{definition}

Let $N$ be a finite von Neumann algebra with a faithful normal trace. The following is directly obtained by combining Lemmas 2.4 and 2.5 in \cite{PoSa}, although these Lemmas were proved for Bernoulli shift actions on standard probability space. The following can be also regarded as a weak $1$-cocycle version of Proposition 3.2 in \cite{Po4}.

\begin{proposition}\label{Weak Coboundary}
Let $G$ be a countable discrete group. Let $\beta$ be a trace preserving weakly mixing action of $G$ on $N$.
A weak $1$-cocycle $\{w_g\}_{g \in G} \subset N$ for $\beta$ is a weak $1$-coboundary if only if there exists a non-zero element $\widetilde{x_0} \in N \otimes N$ satisfying
\begin{eqnarray*}
   (w_g \otimes 1)(\beta_g \otimes \beta_g)(\widetilde{x_0})(1 \otimes w^*_g) = \widetilde{x_0}, \quad g \in G.
\end{eqnarray*}
\end{proposition}

The following is a weak $1$-cocycle version of Proposition 3.6.$3^\circ$ in \cite{Po4}.

\begin{proposition}\label{Extension of Weak Coboundary}
Let $\Gamma$ be a countable discrete group and $G$ be a normal subgroup of $\Gamma$. The group $\Gamma$ acts on a finite von Neumann algebra $N$ in a trace-preserving way by $\beta$. Suppose that the restriction of $\beta$ to $G$ is weakly mixing. Let $\{w_\gamma\}_{\gamma \in \Gamma}$ be a weak $1$-cocycle for $\beta$. If $w|_G$ is a weak $1$-coboundary, then $w$ is a weak $1$-coboundary for the $\Gamma$-action.
\end{proposition}

\begin{proof}
Suppose that $w|_G$ is a weak $1$-coboundary, that is, there exists a unitary element $v$ in $N$ such that $w_g = v \beta_g (v^*) \mod{\mathbb{T}}$ for $g \in G$.
It suffices to show that $\{w^\prime_\gamma\} = \{v^* w_\gamma \beta_\gamma(v)\}$ is in $\mathbb{T}$ for all $\gamma \in \Gamma$. Take arbitrary $\gamma \in \Gamma, g \in G$. Write $h = \gamma^{-1} g \gamma \in G$. Let $\pi_\gamma$ be the unitary on $L^2(N)$ induced from $\beta_\gamma$. Since $w_{h}^\prime, w_g^\prime \in \mathbb{T}$, we get
\begin{eqnarray*}
    w^\prime_\gamma \pi_g {w^\prime_\gamma}^*
    = (w^\prime_\gamma \pi_\gamma)(w^\prime_h \pi_h)(w^\prime_\gamma \pi_\gamma)^*
    = w^\prime_{\gamma h \gamma^{-1}} \pi_{\gamma h \gamma^{-1}}
    = \pi_g
     \ \mod \mathbb{T},
\end{eqnarray*}
By applying these operators to $\hat{1} \in \widehat{N} \subset L^2(N)$, it follows that $w^\prime_\gamma \beta_g({w^\prime_\gamma}^*) \in \mathbb{T}$. Since the $G$-action is weakly mixing, we have $w^\prime_\gamma \in \mathbb{T}$.
\end{proof}

By using the above propositions, we will ``untwist'' some weak $1$-cocycles later. We require some ergodicity assumption on the weak $1$-cocycles.
\begin{definition}
Let $\Gamma$ be a discrete group and $G$ be a subgroup of $\Gamma$. Suppose that its restriction to $G$ is ergodic. Let $\beta$ be a trace preserving action of $\Gamma$ on $N$. A weak $1$-cocycle $w = \{w_g\}_{g \in \Gamma}$ for $\beta$ is said to {\rm be ergodic on $G$}, if the action $\beta^w$ of $G$ is still ergodic, where $\beta^w$ is defined by $\beta^w_g = {\rm Ad} w_g \circ \beta_g, \quad g \in G$.
\end{definition}

Let $\beta$ be a $\Gamma$-action on $N$. Suppose that the diagonal action $\beta \otimes \beta$ on $(N \otimes N, \tr \otimes \tr)$ has an extension $\widetilde{\beta}$ on a finite von Neumann algebra $(\widetilde{N}, \tau)$. The algebra $\widetilde{N}$ is not necessarily identical with $N \otimes N$. When the action $\widetilde{\beta}$ is ergodic on a normal subgroup $G \subset \Gamma$, we get the following:

\begin{proposition}\label{1-cocycle}
Let $\{w_{\gamma}\}_{\gamma \in \Gamma} \subset N$ be a weak $1$-cocycle for $\beta$. Let $\alpha$ be a trace preserving continuous action of $\R$ on $\widetilde{N}$ satisfying the following properties:
\begin{itemize}
   \item
   $\alpha_1(x \otimes 1) = 1 \otimes x$, \quad for all $x \in N$.
   \item
   $\alpha_t \circ \widetilde{\beta}(\gamma)(\widetilde{x}) = \widetilde{\beta}(\gamma) \circ \alpha_t(\widetilde{x})$, \quad for all $t \in \R$, $\gamma \in \Gamma$ and $\widetilde{x} \in \widetilde{N}$.
\end{itemize}
Suppose that the weak $1$-cocycle $\{ w_\gamma \otimes 1 \} \subset \widetilde{N}$ is ergodic for the $G$-action $\widetilde{\beta}|_G$. If the group inclusion $G \subset \Gamma$ has the relative property (T) of Kazhdan, then there exists a non-zero element $\widetilde{x_0} \in \widetilde{N}$ so that $(w_g \otimes 1)\widetilde{\beta}_g(\widetilde{x_0})(1 \otimes w^*_g) = \widetilde{x_0}, \  g \in G$.
\end{proposition}

This is proved in the same way for Bernoulli shift actions on the infinite tensor product of abelian von Neumann algebras (\cite{PoSa}, Lemma 3.5). Since we are interested in actions on the AFD II$_1$ factor, we require the ergodicity assumption on weak $1$-cocycle $\{w_\gamma \otimes 1\}$. For the self-containedness and in order to make it clear where the ergodicity assumption works, we write down a complete proof.

\begin{proof}
For $t \in (0,1]$, let $K_t$ be the convex weak closure of
\begin{eqnarray*}
\{ (w_g \otimes 1)\alpha_t(w_g^* \otimes 1) \ | \ g \in G\} \subset \widetilde{N}
\end{eqnarray*}
and $\widetilde{x_t} \in K_t$ be the unique element whose $2$-norm is minimum in $K_t$. Since
\begin{eqnarray*}
& & (w_g \otimes 1) \widetilde{\beta_g}((w_{g_1} \otimes 1) \alpha_t(w_{g_1}^* \otimes 1)) \alpha_t(w_g^* \otimes 1) \\
&=& (w_g \beta_g(w_{g_1}) \otimes 1) \alpha_t(\beta_g(w_{g_1}^*)w_g^* \otimes 1)\\
&=& (w_{g g_1} \otimes 1) \alpha_t(w_{g g_1}^* \otimes 1), \quad g, g_1 \in G,
\end{eqnarray*}
we have $(w_g \otimes 1)\widetilde{\beta}_g(K_t) \alpha_t(w_g^* \otimes 1) = K_t$, for $g \in G$.
By the uniqueness of $\widetilde{x_t}$, we get
\begin{eqnarray}\label{equation; fixed}
(w_g \otimes 1)\widetilde{\beta}_g(\widetilde{x_t})\alpha_t(w_g^* \otimes 1) = \widetilde{x_t}, \quad g \in G.
\end{eqnarray}
By the assumption, the action $( \mathrm{Ad}(w_g \otimes 1) \circ \widetilde{\beta}_g )_{g \in G}$ is ergodic on $\widetilde{N}$. By the calculation
\begin{eqnarray*}
& & (w_g \otimes 1)\widetilde{\beta}_g ( \widetilde{x_t}\widetilde{x_t}^* ) (w_g^* \otimes 1) \\
&=&  (w_g \otimes 1)\widetilde{\beta}_g ( \widetilde{x_t}) \alpha_t(w_g^* \otimes 1)  \alpha_t(w_g \otimes 1) \widetilde{\beta}_g  (\widetilde{x_t}^* ) (w_g^* \otimes 1) \\
&=& \widetilde{x_t}\widetilde{x_t}^* , \quad g \in G,
\end{eqnarray*}
we get $\widetilde{x_t}\widetilde{x_t}^* \in \C1$. The element $\widetilde{x_t}$ is a scalar multiple of a unitary in $\widetilde{N}$.

We shall next prove that $\widetilde{x_{1/n}}$ is not zero for some positive integer $n$. The pair $( \Gamma, G)$ has the relative property (T) of Kazhdan. By proposition \ref{Proposition; relative (T)}, we can find a positive number $\delta$ and a finite subset $F \subset \Gamma$ satisfying the following condition: If a unitary representation $(\pi, \mathcal{H})$ of $\Gamma $ and a unit vector $\xi$ of $\mathcal{H}$ satisfy $\| \pi(\gamma) \xi - \xi \| \le \delta \ (\gamma \in F)$, then $\| \pi(g) \xi - \xi \| \le 1/2 \ (g \in G)$. By the continuity of the action $\alpha$, there exists $n$ such that
\begin{eqnarray*}
   \| (w_\gamma \otimes 1) \alpha_{1/n} (w_\gamma^* \otimes 1) - 1 \|_{\tr,2} \le \delta, \quad \gamma \in F,
\end{eqnarray*}
The actions $\beta$ and $(\alpha_{1/n}^l)_{l \in \Z}$ on $\widetilde{N}$ give a $\Gamma \times \Z$ action on $\widetilde{N}$. Let $P$ be the crossed product von Neumann algebra
$P = \widetilde{N} \rtimes (\Gamma \times \Z)$.
Let $(U_\gamma)_{\gamma \in \Gamma}$ and $W$ be the implementing unitaries in $P$ for $\Gamma$ and $1 \in \Z$ respectively. We put $V_{\gamma} = (w_{\gamma} \otimes 1)U_{\gamma}, \  \gamma \in \Gamma$. We regard $\textrm{Ad}V_\cdot$ as a unitary representation of $\Gamma$ on $L^2(P)$. Since
\begin{eqnarray*}
     \| \textrm{Ad}V_{\gamma}(W) - W \|_{L^2(P)}
 &=& \| (w_{\gamma} \otimes 1)W(w_{\gamma}^* \otimes 1)W^* - 1  \|_{L^2(P)}  \\
 &=& \| (w_{\gamma} \otimes 1)\alpha_{1/n}(w_{\gamma}^* \otimes 1) - 1  \|_{L^2(\widetilde{N})}
 \le \delta, \quad \gamma \in F,
\end{eqnarray*}
we have the following inequality:
\begin{eqnarray*}
1 / 2 \ge  \| \textrm{Ad}V_g(W) - W \|_{L^2(P)}
        =    \| (w_g \otimes 1)\alpha_{1/n}(w_g^* \otimes 1) - 1  \|_{L^2(\widetilde{N})}, \quad g \in G.
\end{eqnarray*}
We get $1 / 2 \ge  \| \widetilde{x_{1/n}} - 1  \|_{L^2(\widetilde{N})}$
and $\widetilde{x_{1/n}} \neq 0$.

Let $\widetilde{u_{1/n}}$ be the unitary of $\widetilde{N}$ given by a scalar multiple of $\widetilde{x_{1/n}}$. By equation (\ref{equation; fixed}), the unitary satisfies
\begin{eqnarray*}
(w_g \otimes 1)\widetilde{\beta}_g(\widetilde{u_{1/n}}) \alpha_{1/n}(w_g^* \otimes 1) &=& \widetilde{u_{1/n}}, \quad g \in G.
\end{eqnarray*}
Let $\widetilde{x_0}$ be the unitary defined by
\begin{eqnarray*}
     \widetilde{x_0} = \widetilde{u_{1/n}} \alpha_{1/n}( \widetilde{u_{1/n}})  \alpha_{2/n}( \widetilde{u_{1/n}} ) \hdots \alpha_{(n-1)/n}( \widetilde{u_{1/n}}).
\end{eqnarray*}
By direct computations, we have the following desired equality:
\begin{eqnarray*}
    (w_g \otimes 1)\widetilde{\beta}_g(\widetilde{x_0})(1 \otimes w_g^*)
= (w_g \otimes 1)\widetilde{\beta}_g(\widetilde{x_0}) \alpha_1 (w_g^* \otimes 1)
= \widetilde{x_0}, \quad g \in G.
\end{eqnarray*}
\end{proof}

\begin{theorem}\label{CocycleConjugacyAndConjugacy}
Let $\beta = \beta(H, \mu, \chi)$ be a twisted Bernoulli shift action on $N(H, \mu)$. Suppose that $N(H, \mu)$ is the AFD factor of type $\mathrm{II}_1$ and that there exists a continuous $\R$-action $(\alpha^{(0)}_t)_{t \in \R}$ on $L_\mu(H) \otimes L_\mu(H)$ satisfying the following properties:
\begin{itemize}
    \item
    For any $x \in L_\mu(H)$, $\alpha^{(0)}_1(x \otimes 1) = 1 \otimes x$,
    \item
    The automorphism $\alpha^{(0)}_t$ commutes with the diagonal action of $\widehat{H}$.
\end{itemize}
Let $\beta^{(1)}$ be another action of $\G$ on the AFD factor $N^{(1)}$ of type $\mathrm{II}_1$ and suppose that its restriction to $\Z^2$ is ergodic. The action $\beta^{(1)}$ is outer conjugate to $\beta$, if and only if $\beta^{(1)}$ is conjugate to $\beta$.
\end{theorem}
\begin{proof}
We deduce from outer conjugacy to conjugacy in the above situation. Let $\theta$ be a $*$-isomorphism from $N^{(1)}$ onto $N(H, \mu)$ which gives the outer conjugacy of the action $\beta^{(1)}$ and $\beta = \beta(H, \mu, \chi)$. There exists a weak $1$-cocycle $\{ w_\gamma \}_{\gamma \in \G}$ for $\beta$ satisfying
\begin{eqnarray*}
    \theta \circ \beta^{(1)}(\gamma) = \mathrm{Ad} w_\gamma \circ \beta(\gamma) \circ \theta, \quad \gamma \in \G.
\end{eqnarray*}
Since the action $\beta^{(1)}$ is ergodic on $\Z^2$, the weak $1$-cocycle $w$ is ergodic on $\Z^2$.

We use the notations $\Gamma_0, G_0$ given in Section \ref{definition}. Let $\widetilde{\rho}$ be the diagonal action $\rho \otimes \rho$ of $\Gamma_0$ on the tensor product algebra $\widetilde{M} = L_{\widetilde\mu}(\oplus_{\Z^2} H) \otimes L_{\widetilde\mu}(\oplus_{\Z^2} H)$:
\begin{eqnarray*}
    \widetilde{\rho}(\gamma_0) (a \otimes b) = \rho(\gamma_0)(a) \otimes \rho(\gamma_0)(b).
\end{eqnarray*}
The fixed point algebra $\widetilde{N} \subset \widetilde{M}$ of the diagonal $\widehat{H}$-action contains $N(H, \mu) \otimes N(H, \mu)$. Since $\G = \Gamma_0 / \widehat{H}$, the action $\widetilde\rho$ gives a $\G$-action $\widetilde\beta$ on $\widetilde{N}$. The action $\widetilde\beta$ is the extension of the diagonal action $\beta \otimes \beta$ on $N(H, \mu) \otimes N(H, \mu)$. We denote by $\alpha_t$ the action on $\widetilde{M} \cong \bigotimes_{\Z^2} (L_\mu(H) \otimes L_\mu(H))$ given by the infinite tensor product of the $\R$-action $\alpha_t^{(0)}$.
By the assumption on $\alpha^{(0)}_t$, the $\R$-action $\alpha_t$ commutes with the action $\widetilde{\rho}$. It follows that the subalgebra $\widetilde{N}$ is globally invariant under $\alpha_t$.

The set of unitary $\{W_\gamma = w_\gamma \otimes 1\}_{\gamma \in \G} \subset \widetilde{N}$ is a weak $1$-cocycle for $\widetilde\beta$. We shall prove that this weak $1$-cocycle is ergodic on $\Z^2$.
Let $a$ be an element in $\widetilde{N}$ fixed under $\widetilde\beta |_{\Z^2}$. The element $a$ can be written as
$a = \sum_{\lambda \in \oplus_{\Z^2} H} a_\lambda \otimes u(\lambda)$
in $L^2 \widetilde{M}$, where $a_\lambda \otimes 1 = E_{M \otimes \C} (a (1 \otimes u(\lambda))^*)$. Since $a$ is fixed under the action of $\Z^2$, we have
\begin{eqnarray*}
    a &=& \widetilde\beta^{W}(k)(a)
      = \sum_{\lambda \in \oplus_{\Z^2} H} \mathrm{Ad}w_k \circ \rho(1, k) (a_\lambda) \otimes \rho(1,k) (u(\lambda))\\
      &=& \sum_{\lambda \in \oplus_{\Z^2} H}
      \mathrm{Ad}w_k \circ \rho(1, k) (a_\lambda) \otimes
      \prod_{l \in \Z^2} \chi(\lambda(l))^{\det(k, l)} u(k \cdot \lambda).
\end{eqnarray*}
Since $\mathrm{Ad}w_k \circ \rho(1, k)$ preserves the $2$-norm, we get $\| a_\lambda \|_2 = \| a_{k^{-1} \cdot \lambda} \|_2$. Since $\| a \|_2^2 = \sum \| a_\lambda \|_2^2 < \infty$ and the set $\{ a_{k^{-1} \cdot \lambda} \ | \ k \in \Z^2\}$ is infinite for $\lambda \neq 0$, it turns out that $a_\lambda = 0$ for $\lambda \neq 0$ and thus $a \in \widetilde{N} \cap (M \otimes \C) = N(H, \mu) \otimes \C$. By the ergodicity of the $\Z^2$-action $\{\mathrm{Ad}w_k \circ \beta(k)\}$, we get $a \in \C$. We conclude that the weak $1$-cocycle $\{W_\gamma\} \subset \widetilde{N}$ is ergodic on $\Z^2$.

By the relative property (T) for the inclusion $\Z^2 \subset \G$ and Proposition \ref{1-cocycle}, there exists a non-zero element $\widetilde{x_0} \in \widetilde{N}$ satisfying
\begin{eqnarray*}
    (w_k \otimes 1)\widetilde{\beta}(k)(\widetilde{x_0})(1 \otimes w^*_k) = \widetilde{x_0}, \  k \in \Z^2.
\end{eqnarray*}
The element $\widetilde{x_0}$ can be written as the following Fourier expansion:
\begin{eqnarray*}
    \widetilde{x_0} = \sum c(\lambda_1, \lambda_2) u(\lambda_1) \otimes u(\lambda_2) \in \widetilde{N} \subset L_{\widetilde\mu}(\oplus_{\Z^2}H) \otimes L_{\widetilde\mu}(\oplus_{\Z^2}H).
\end{eqnarray*}
Here $c(\lambda_1, \lambda_2)$ is a complex number and $(\lambda_1, \lambda_2) \in (\oplus_{\Z^2} H)^2$ runs through all pairs satisfying $\sum_{k \in \Z^2} (\lambda_1(k) + \lambda_2(k)) = 0$. Choose and fix a pair $(\lambda_1, \lambda_2)$ satisfying
\begin{eqnarray*}
- \sum_{k \in \Z^2} \lambda_1(k) = h = \sum_{k \in \Z^2} \lambda_2(k), \quad c(\lambda_1, \lambda_2) \neq 0.
\end{eqnarray*}
Let $v_h^\prime \in M$ be the unitary written as $ v_h^\prime = u_h \otimes 1 \otimes 1 \otimes \cdots$, where $u_h \in L_{\mu}(H)$ is placed on $0 \in \Z^2$. The following unitaries $\{w^\prime_\gamma\} \subset N(H, \mu)$ give a weak $1$-cocycle for $\beta$:
\begin{eqnarray*}
    w^\prime_{(k, \gamma_0)} = v_h^\prime w_{(k, \gamma_0)} \rho(1, k, \gamma_0) ({v_h^\prime}^*), \quad (k, \gamma_0) \in \G.
\end{eqnarray*}
Letting $\widetilde{y} = (v_h^\prime \otimes 1) \widetilde{x_0} (1 \otimes v_h^\prime)^* \in \widetilde{M}$, we get
\begin{eqnarray*}
    \widetilde{y}    &=& (w^\prime_k \otimes 1)\widetilde{\beta}(k)(\widetilde{y})(1 \otimes {w_k^\prime}^*), \quad k \in \Z^2.
\end{eqnarray*}
Applying the trace preserving conditional expectation $E = E_{N(H, \mu) \otimes N(H, \mu)}$, we get
\begin{eqnarray*}
    E(\widetilde{y}) &=& (w^\prime_k \otimes 1)E(\widetilde{\beta}(k)(\widetilde{y}))(1 \otimes {w_k^\prime}^*)\\
         &=& (w^\prime_k \otimes 1)\widetilde{\beta}(k)(E(\widetilde{y}))(1 \otimes {w_k^\prime}^*), \quad k \in \Z^2.
\end{eqnarray*}
Since the Fourier coefficient of $\widetilde{x_0}$ at $(\lambda_1, \lambda_2) \in (\oplus_{\Z^2} H)^2$ is not zero, that of $E(\widetilde{y})$ at $(\lambda_1 + \delta_{h, 0}, \lambda_2 - \delta_{h, 0}) \in \Lambda(H)^2$ is also non-zero, where $\delta_{h, 0} \in \oplus_{\Z^2} H$ is zero on $\Z^2 \setminus \{0\}$ and is $h$ on $0 \in \Z^2$.
By Proposition \ref{Weak Coboundary}, it follows that the weak $1$-cocycle $\{w^\prime_{(k, e)} \}_{k \in \Z^2} \subset N(H, \mu)$ is a weak $1$-coboundary of $\beta |_{\Z^2}$. Since the $\Z^2$-action $\beta |_{\Z^2}$ is weakly mixing, $w^\prime$ is a weak $1$-coboundary on $\G$, by Proposition \ref{Extension of Weak Coboundary}. In other words, there exists $v \in N(H, \mu)$ satisfying
\begin{eqnarray*}
    w^\prime_\gamma &=& v \beta(\gamma)(v^*) \ \mod \mathbb{T},\\
    w_\gamma &=& {v_h^\prime}^* v \rho(1, \gamma)(v^* v_h^\prime) \ \mod \mathbb{T}, \quad \gamma \in \G.
\end{eqnarray*}
Noting that $u = v^* v_h^\prime \in M$ is a normalizer of $N(H, \mu)$, we get
\begin{eqnarray*}
    (\mathrm{Ad}(u) \circ \theta) \circ \beta^{(1)}(\gamma)
    &=& \mathrm{Ad}(u) \circ \mathrm{Ad}(w_\gamma) \circ \beta(\gamma) \circ \theta\\
    &=& \mathrm{Ad}(\rho(1, \gamma)(u)) \circ \beta(\gamma) \circ \theta\\
    &=& \rho(1, \gamma) \circ \mathrm{Ad}(u) \circ \theta\\
    &=& \beta(\gamma) \circ (\mathrm{Ad}(u) \circ \theta), \quad \gamma \in \G.
\end{eqnarray*}
Thus we get the conjugacy of two $\G$-actions $\beta^{(0)}$ and $\beta$.
\end{proof}

We can always apply Theorem \ref{CocycleConjugacyAndConjugacy} if $H$ is finite.

\begin{corollary}\label{Finite Abelian}
Let $H$ be a finite abelian group and let $\beta = \beta(H, \mu, \chi)$ be a twisted Bernoulli shift action on $N(H, \mu)$. Suppose that $N(H, \mu)$ is the AFD factor of type $\mathrm{II}_1$. Let $\beta^{(1)}$ be an action of $\G$ on the AFD factor $N^{(1)}$ of type $\mathrm{II}_1$ and suppose that its restriction to $\Z^2$ is ergodic. The action $\beta^{(1)}$ is outer  conjugate to $\beta$, if and only if $\beta^{(1)}$ is conjugate to $\beta$.
\end{corollary}

\begin{proof}
We have only to construct an $\R$-action on $L_\mu(H) \otimes L_\mu(H)$ satisfying the properties in Theorem \ref{CocycleConjugacyAndConjugacy}. Let $U$ be an element of $L_\mu(H) \otimes L_\mu(H)$ defined by
\begin{eqnarray*}
    U = \frac{1}{|H|^{1/2}}\sum_{h \in H} u_h \otimes u_h^*.
\end{eqnarray*}
We note  that $\mu^*\mu(g, \cdot)$ is a character of $H$ and that it is not identically $1$ provided $g \neq 0$ by Lemma \ref{AFD II_1}. The element $U$ is self-adjoint and unitary, since
\begin{eqnarray*}
    U^* &=& \frac{1}{|H|^{1/2}}\sum_{h \in H} u_h^* \otimes u_h
           = \frac{1}{|H|^{1/2}}\sum_{h \in H} \overline{\mu(h, -h)}u_{-h} \otimes \mu(h, -h) u_{-h}^* = U,\\
    U^2 &=& \frac{1}{|H|}\sum_{g, h \in H} u_g u_h\otimes u_g^* u_h^*
           =  \frac{1}{|H|}\sum_{g, h \in H} \mu^*\mu(g, h) u_{g + h} \otimes u_{g + h}^*\\
          &=& \frac{1}{|H|}\sum_{g \in H} \left(\sum_{h \in H} \mu^*\mu(g, h - g)\right) u_g \otimes u_g^* = 1.
\end{eqnarray*}
The operator $U$ is a fixed point under the action of $\widehat{H}$, so the projections $P_1 = (1 + U)/2 $ and $P_{-1} = (1 - U) / 2$ are also fixed points. Thus the $\R$-action $\alpha^{(0)}_t = \mathrm{Ad}(P_1 + \exp{(i \pi t)} P_{-1})$ commutes with the $\widehat{H}$-action.
The automorphism $\alpha^{(0)}_1$ satisfies
\begin{eqnarray*}
    \alpha^{(0)}_1 (u_g \otimes 1) = U (u_g \otimes 1) U^* = (1 \otimes u_g) U U^* = 1 \otimes u_g, \quad g \in H.
\end{eqnarray*}
This verifies the first condition for $\alpha^{(0)}$.
\end{proof}

\begin{corollary}
Let $Q$ be a set consisting of odd prime numbers and $\beta_Q$ be the twisted Bernoulli shift action defined in Corollary \ref{Collection of Odd Prime numbers}. Let $\beta$ be a $\G$-action on the AFD factor of type $\mathrm{II}_1$ whose restriction to $\Z^2$ is ergodic. The actions $\beta_Q$ and $\beta$ are outer conjugate if and only if they are conjugate. In particular, $\{ \beta_Q \}$ is an uncountable family of $\G$-actions which are mutually non outer conjugate.
\end{corollary}
\begin{proof}
We will use the notation given in Corollary \ref{Odd Prime numbers} and \ref{Collection of Odd Prime numbers}. Let $\alpha^{(q)}_t$ be the $\R$-action on $L_{\mu_q}(H_q) \otimes L_{\mu_q}(H_q)$ constructed as in the previous corollary. We define the $\R$-action $\alpha^{(Q)}$ on $L_{\mu_Q}(H_Q) \otimes L_{\mu_Q}(H_Q)$ by $\alpha^{(Q)}_t(\otimes_{q \in Q} x_q) = \otimes_{q \in Q} \alpha^{(q)}_t(x_q)$, where $x_q \in L_{\mu_q}(H_q) \otimes L_{\mu_q}(H_q)$ and $x_q \neq 1$ only for finitely many $q$. The $\R$-action satisfies the conditions in Theorem \ref{CocycleConjugacyAndConjugacy}.
By Corollary \ref{Collection of Odd Prime numbers}, $\{ \beta_Q \}$ are mutually non conjugate and their restriction to $\Z^2$ is ergodic. Thus they are mutually non outer conjugate.
\end{proof}

\begin{acknowledgment}
The author would like to thank Professor Yasuyuki Kawahigashi for helpful conversations. He thanks the referee for careful reading and numerous detailed comments. He is supported by JSPS Research Fellowships for Young Scientists.
\end{acknowledgment}

\end{document}